\documentclass[12pt,a4paper]{article}
\usepackage[T1]{fontenc}
\usepackage{multirow,scalerel}
\usepackage{setspace}

\usepackage[numbered]{bookmark}

\usepackage{appendix}

\usepackage{amsfonts,enumerate,amsmath}
\usepackage{latexsym}
\usepackage{color}
\usepackage{natbib}
\usepackage{array}

\usepackage{verbatim}
\usepackage{graphicx}
\usepackage{array}
\usepackage{float,bm}
\usepackage{multirow}

\renewcommand{\theequation}{\arabic{section}.\arabic{equation}}

\newcommand{\qed}{\hfill\Box}

\graphicspath{{./images/}}

\renewcommand{\cite}{\citet*}

\newcommand{\bqn}{\begin{eqnarray*}}
\newcommand{\eqn}{\end{eqnarray*}}

\DeclareMathOperator{\tr}{tr}         
\DeclareMathOperator{\im}{Im}      \DeclareMathOperator{\diag}{diag}

\numberwithin{equation}{section}  

\newtheorem{thm}{Theorem}

\newtheorem{rem}{Remark}
\newtheorem{cor}{Corollary}

\newtheorem{assp}{Assumption}

\hypersetup{
  colorlinks   = true, 
  urlcolor     = blue, 
  linkcolor    = blue, 
  citecolor   = blue 
}

\begin{document}

\title{
High-dimensional covariance matrices under dynamic volatility models: asymptotics and shrinkage estimation
}

\author{Yi Ding \thanks{Faculty of Business Administration, University of Macau, Macau. Email: yiding@um.edu.mo}
 \and Xinghua Zheng \thanks{Department of ISOM, Hong Kong University of Science and Technology, Clear Water Bay, Kowloon, Hong Kong. Email: xhzheng@ust.hk}}
\date{November 18, 2022}

\maketitle
\begin{abstract}
\noindent

We study the estimation of the high-dimensional covariance matrix and its eigenvalues under dynamic volatility models. Data under such models have nonlinear dependency both cross-sectionally and temporally. We first investigate the empirical spectral distribution (ESD) of the sample covariance matrix under scalar~BEKK models and  establish conditions under which the limiting spectral distribution (LSD) is either the same as  or different from the~\mbox{i.i.d.} case. We then propose a time-variation adjusted (TV-adj) sample covariance matrix and prove that its LSD follows the same Mar$\check{\text{c}}$enko-Pastur law as the~\mbox{i.i.d.} case. Based on the~asymptotics of the~TV-adj sample covariance matrix, we develop a consistent population spectrum estimator and an asymptotically optimal~nonlinear shrinkage estimator of the unconditional covariance matrix.  

 \vskip 0.3cm



\end{abstract}

\section{Introduction}
\subsection{The Mar$\check{\text{c}}$enko-Pastur law for the sample covariance matrix}
Random matrix theory (RMT) is a powerful tool in the study of high-dimensional statistics. When the dimension and the sample size grow  proportionally, for \mbox{i.i.d.} data, it is well known that the limiting spectral distribution (LSD) of the sample covariance matrix is connected to that of the population covariance matrix through the Mar$\check{\text{c}}$enko-Pastur equation; see, e.g.,  \cite{MP67}, \cite{yin1986limiting}, \cite{SB95}, and \cite{silverstein1995strong}. 
Algorithms of recovering population spectrum based on the Mar$\check{\text{c}}$enko-Pastur (M-P) law have been developed in \cite{el2008spectrum},  \cite{LW12,LW15,LW20}. All these studies  focus on the case where observations are \mbox{i.i.d.}. 
\subsection{Dynamic volatility models}
An important feature of  financial returns is that their volatilities are time-varying and dependent over time. Dynamic volatility models such as the multivariate~GARCH [\cite{engle1984combining}, \cite{bollerslev1988capital}], the~DCC model [\cite{engle2002dynamic}] and the~BEKK model [\cite{engle1995multivariate}] are popular in studying the dynamic variances and covariances. In particular, the widely used scalar~BEKK model describes the dynamics of covariance matrix as follows:
\begin{equation}\label{BEKK}
\boldsymbol{\Sigma}_{t+1}=(1-a-b)\overline{\boldsymbol{\Sigma}}+a\mathbf{R}_{t}\mathbf{R}_{t}^T+b\boldsymbol{\Sigma}_{t},
\end{equation}
where~$\overline{\boldsymbol{\Sigma}}$ is the unconditional covariance matrix,~$\boldsymbol{\Sigma}_{t}$ is the conditional covariance matrix of the returns~$\mathbf{R}_t=(R_{1t}, ..., R_{pt})^T$, and~$0\leq a,b\leq 1$ with $a+b<1$  are related parameters. The parameter $a$ is sometimes referred to as the innovation coefficient, and $a+b$ the persistence coefficient.

To estimate a dynamic volatility model, a common approach is  variance/correlation targeting [\cite{francq2011merits}, \cite{pedersen2014multivariate}, \cite{pakel2020fitting}]. The method requires estimating the unconditional covariance/correlation matrix. When the dimension is high,  the sample covariance/correlation matrix works poorly. For large dynamic volatility models, estimating the unconditional covariance/correlation matrix is  challenging and calls for rigorous investigation.

\subsection{Existing research in RMT for sample covariance matrix when there is time dependency}
There exists a growing literature on the study of limiting spectral properties of the sample covariance matrix when there is time  dependency. \cite{jin2009limiting}, \cite{yao2012note}, \cite{liu2015marvcenko}, and \cite{BB16} obtain the~LSD of the sample covariance/autocovariance matrix  of linearly dependent time series that can be transformed into data with independent columns. 
\cite{banna2015limiting} and \cite{merlevede2016empirical} investigate the LSD of the sample covariance matrix of stationary dependent processes with independent rows. \cite{yaskov2018spectrum} focuses on the case where data have  dependence in finite lags.
\cite{ZL11} establish the~LSD, and \cite{YZC20} derive the central limit theorem of linear spectral statistics of sample covariance matrix under elliptical models with time-varying co-volatilities.

 \cite{ELW19} propose to estimate the unconditional covariance/correlation matrix under large~BEKK/DCC models using the nonlinear shrinkage (NLS) estimator developed in \cite{LW12,LW15}. It has been documented in  \cite{LW12,LW15} that the NLS estimator has several advantages in estimating the high-dimensional covariance matrix. For example, it is structure-free, and more importantly, for \mbox{i.i.d.} data,  it is consistent in estimating the asymptotically optimal shrinkage estimator in the class of rotation-equivariant estimators; see \cite{LW12,LW15} for detailed explanations. 
It is worth emphasizing that the asymptotic property of the NLS relies on that the LSD of the 
sample covariance matrix follows the M-P law. For large dynamic volatility models,  whether the~NLS still enjoys the desirable asymptotic property is unclear.

 \subsection{Our contributions}
 
 We aim to estimate the unconditional covariance matrix under large dynamic volatility models. An important and natural question motivated by the proposal of \cite{ELW19} is:  
 \emph{does the NLS work under large dynamic volatility models?}

\begin{figure}[H]
\centering
\includegraphics[width=0.4\textwidth]{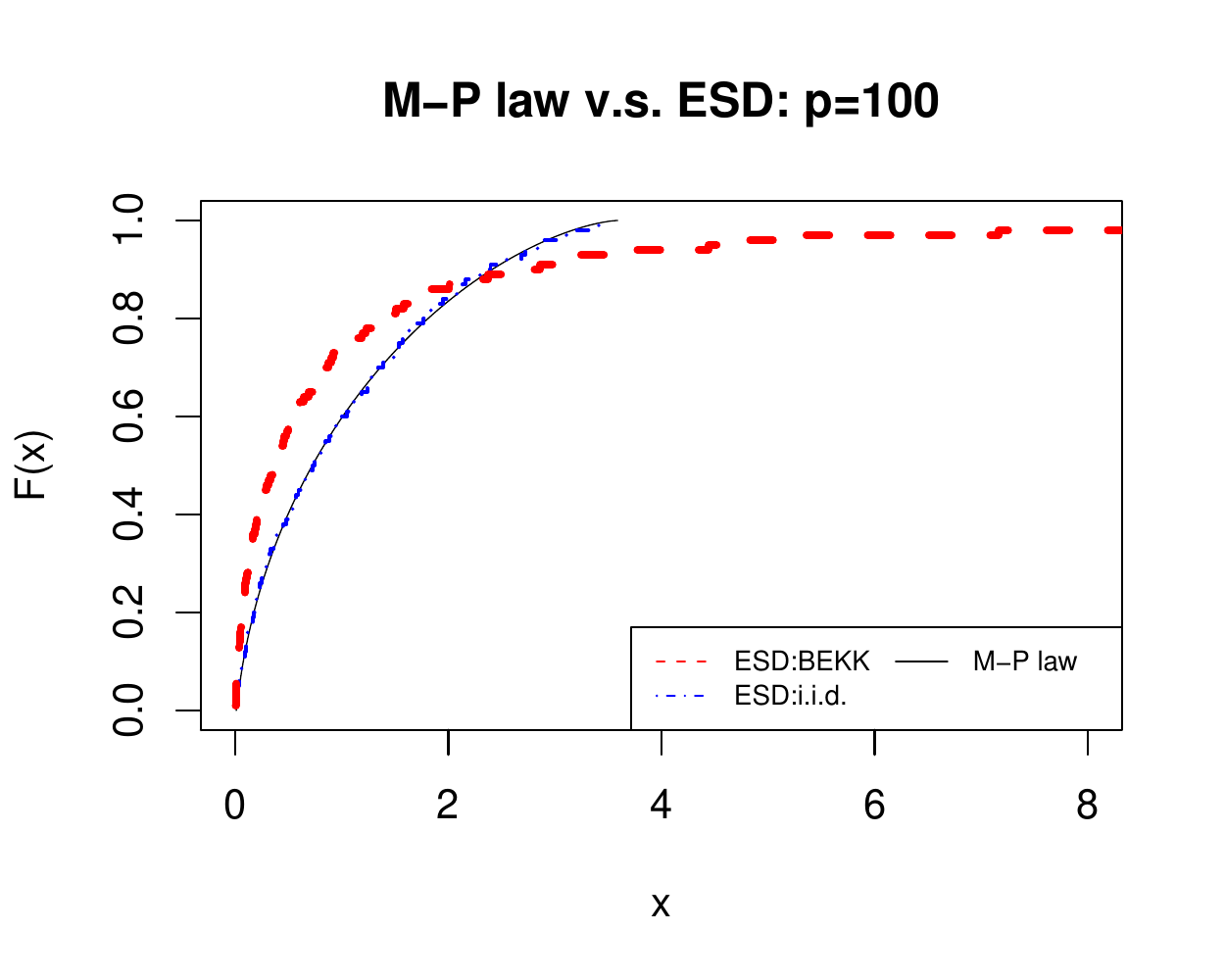}
\includegraphics[width=0.4\textwidth]{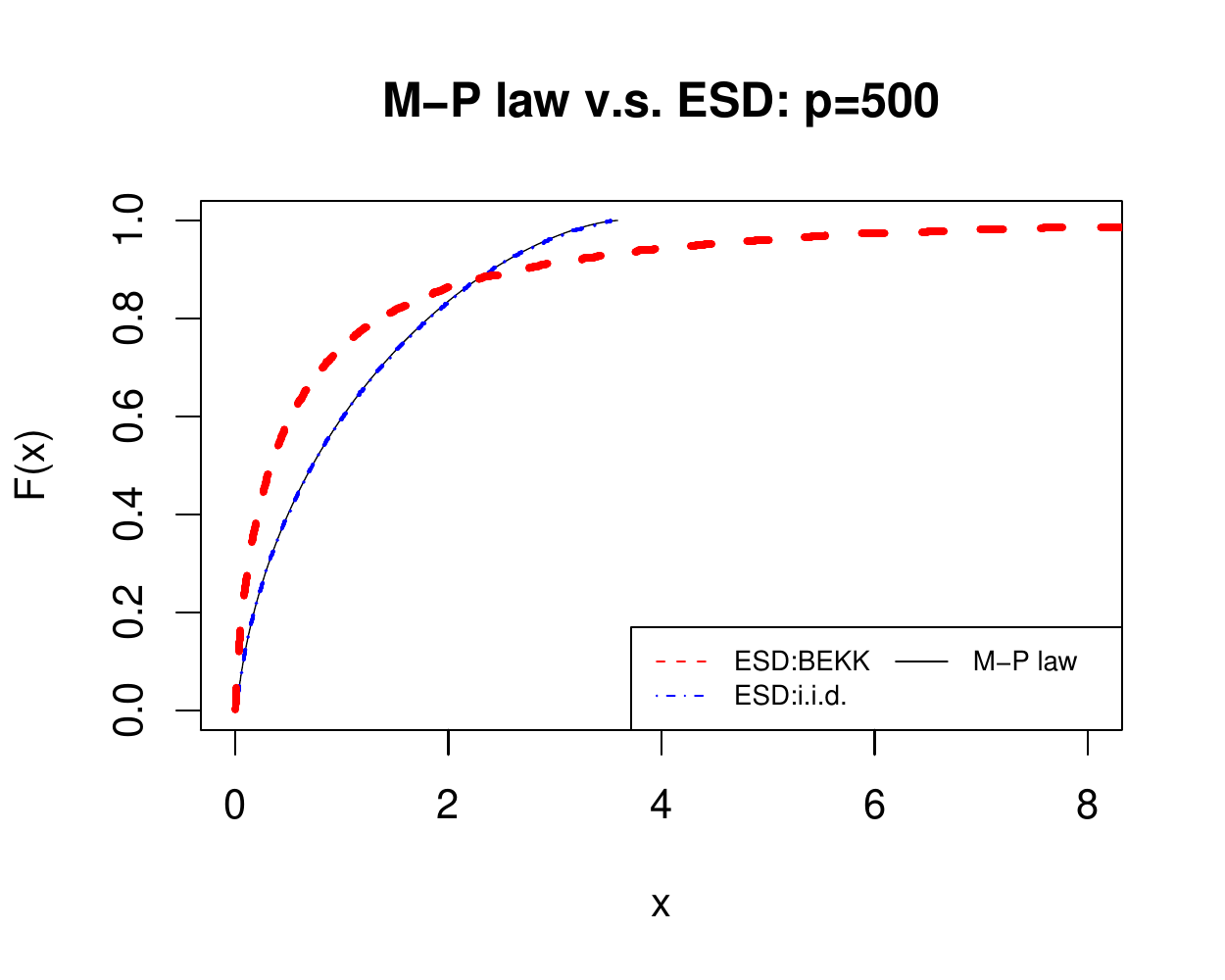}
\caption{{Empirical spectral distributions of sample covariance matrices under~the BEKK model, in comparison with the ESD based on \mbox{i.i.d.} data and the Mar$\check{\text{c}}$enko-Pastur (M-P) distribution. }}\label{Fig1}
\end{figure}

To see how the dynamic volatility model can affect the spectral distribution of the sample covariance matrix,
we simulate data from BEKK model~\eqref{BEKK} with \mbox{$\overline{\boldsymbol{\Sigma}}=\mathbf{I}$},~$a=0.05$ and~\mbox{$b=0.9$}, which is the setting used in \cite{ELW19}. The dimensions~$p=100, 500$, and the sample size~$n$ satisfies~$p/n=0.8$. We compute the empirical spectral distribution (ESD) of  the sample covariance matrix and compare it with the Mar$\check{\text{c}}$enko-Pastur (M-P) distribution. 
The results are shown in~Figure~\ref{Fig1}.
\noindent Figure~\ref{Fig1} shows that the ESD of the sample covariance matrix under the~BEKK model behaves differently from the \mbox{M-P} law. Therefore, it is problematic to perform NLS on the sample covariance matrix the same way as in the \mbox{i.i.d.} case.

In this paper, we investigate the limiting spectral properties of the sample covariance matrix under large~BEKK models. We show that if {$\eta(a,b,p):=\big(a/(1-a-b)\big)\min\big(\sqrt{p(1-a-b)},1\big)\to 0$}, then the~LSD of the sample covariance matrix shares the same limit as the \mbox{i.i.d.} case; see~Theorem~\ref{reducable}. We call this case the reducible case. On the other hand, if~${\eta(a,b,p)}$ is bounded away from zero, then the spectral distribution of the sample covariance under the~BEKK model is more heavy-tailed than the \mbox{i.i.d.} case; see~Theorem~\ref{nonreducable}.  We call this case the non-reducible case.

Next, we address the problem of population spectrum estimation under large~BEKK models. We first estimate the  parameters~$a$ and~$b$ using a QMLE from univariate~GARCH models. Next, we develop a projection matrix,~$\mathbf{P}_t$, which can track the time-variation in the conditional covariance matrices. We then
define time-variation adjusted returns,~$\widetilde{\mathbf{R}}_t=\mathbf{P}_t^{-1/2}\mathbf{R}_t$, and a \emph{time-variation adjusted  (TV-adj) sample covariance matrix}~$\widetilde{\mathbf{S}}_n=\sum_{t=1}^n\widetilde{\mathbf{R}}_t\widetilde{\mathbf{R}}_t^T/n$. We prove that
the~TV-adj sample covariance matrix shares the same~LSD as the \mbox{i.i.d.} case; see Theorem~\ref{recovered}. Using the~TV-adj sample covariance matrix and existing Mar$\check{\text{c}}$enko-Pastur law reversing algorithms, we obtain a~TV-adj shrinkage estimator of the population spectrum, which is show to be consistent; see Corollary \ref{NLS_ADJ} for the exact statement.  

Finally, we tackle the problem of unconditional covariance matrix estimation under large~BEKK models. We develop a TV-adj nonlinear shrinkage (NLS) estimator and show that it consistently estimates the asymptotically optimal shrinkage estimator; see Theorem \ref{NLS_ADJ_cov}.
\vskip 0.15cm
In summary, our contributions lie in the following aspects:
\begin{itemize}
\item First, we establish conditions under which the LSD of the sample covariance matrix under large~BEKK models is the same as or different from the \mbox{i.i.d.} case. 
\item Second, we propose a~TV-adj sample covariance matrix and develop an estimator that can consistently recover the population eigenvalues under large~BEKK models. 
\item Third, we develop a~\mbox{TV-adj NLS} estimator and prove that it is asymptotically optimal.  
\end{itemize}

 The rest of this paper is organized as follows. The main theoretical results are given in Section \ref{main}. Simulation studies are presented in Section \ref{Simu}.  We conclude in Section \ref{Conc}. The proof of Theorem \ref{recovered} is presented in Section \ref{proof_recovered}. The proofs of other main results are collected in the Supplementary Material \cite{DZ22_supp}.
 
 The following notation is used throughout the paper.  
For any matrix~$\mathbf{A}=(A_{ij})$, its spectral norm is defined as~$\|\mathbf{A}\|=\max_{\|\mathbf{x}\|\leq 1}\sqrt{\mathbf{x}^T\mathbf{A}^T\mathbf{A}\mathbf{x}}$, where $\|\mathbf{x}\|=\sqrt{\sum x_i^2}$ for any vector $\mathbf{x}=(x_i)$; the Frobenius norm is defined as~$\|\mathbf{A}\|_{F}=\sqrt{\sum_{i,j}A_{ij}^2}$.  We write $\mathbf{A}\geq 0 (>0)$ if the matrix $\mathbf{A}$ is positive semi-definite (positive definite), and $\mathbf{A}\geq \mathbf{B} (>\mathbf{B})$ if $\mathbf{A}-\mathbf{B}\geq 0 (>0)$. If $\mathbf{A}\geq 0$, $\mathbf{A}^{1/2}$ is defined as the matrix that satisfies $\mathbf{A}^{1/2}\geq 0$ and \mbox{$(\mathbf{A}^{1/2})^2=\mathbf{A}$}. For any symmetric matrix $\mathbf{A}$ with eigenvalues $\lambda_1, ..., \lambda_p$, its empirical spectral distribution (ESD)  is defined as $F^{\mathbf{A}}(x)=\sum_{j=1}^p\mathbf{1}_{[\lambda_j,+\infty)}(x)/p,$   $x\in \mathbb{R}$. Finally, we use $\overset{\text{P}}\to$ to represent convergence in probability.
 
\section{Main Results}\label{main}
\subsection{Setting and assumptions}
Under a dynamic volatility model, returns are modeled as~$\mathbf{R}_t=(\boldsymbol{\Sigma}_{t})^{1/2}\mathbf{z}_t$, where~$\mathbf{z}_t=(z_{1t}, ..., z_{pt})^T$ are \mbox{i.i.d.} with mean zero and covariance matrix~$\mathbf{I}$. 
We suppose that~$(\mathbf{R}_t)$ follows the scalar~BEKK model \eqref{BEKK}. Define $\mathbf{R}^0_{t}=(\overline{\boldsymbol{\Sigma}})^{1/2}\mathbf{z}_t$, $t=1, ..., T$, which share the same unconditional covariance matrix as $\mathbf{R}_t$ but are \mbox{i.i.d.}. Denote the corresponding sample covariance matrices by 
$$\mathbf{S}_n=\frac{1}{n}\sum_{t=1}^n \mathbf{R}_{t}\mathbf{R}_{t}^T,\,\text{  and }\;\mathbf{S}_n^0=\frac{1}{n}\sum_{t=1}^n \mathbf{R}^0_{t}(\mathbf{R}^0_{t})^T.$$
We write~$\widehat{\lambda}_1\geq ... \geq\widehat{\lambda}_{p}$ as the eigenvalues of~$\mathbf{S}_n$, and~$\widehat{\lambda}^0_1\geq ... \geq\widehat{\lambda}^0_{p}$ as the eigenvalues of~$\mathbf{S}^0_n$. 

We impose the following assumptions.
\begin{assp}\label{asump1}\hfill

\vskip 0.2cm

(i)~{$\mathbf{z}_{t}\underset{\text{i.i.d.}}\sim {N}(0, \mathbf{I})$.}



(ii)~$\overline{\boldsymbol{\Sigma}}$ is nonnegative definite and its~ESD,~$F^{\overline{\boldsymbol{\Sigma}}}$, converges in distribution to a probability distribution~$H$ on $[0, \infty)$ as $p \to \infty$, and $H\neq \delta(0)$, the Dirac measure at~$0$.

(iii)~$\|\overline{\boldsymbol{\Sigma}}\|<C$ for some constant~$C>0$. 

(iv)~The dimension $p$ and the sample size $n$ satisfy that~$p, n\to\infty$, and~\mbox{$p/n\to y>0$}. 

\end{assp}

About the parameters $a$ and $b$, we allow them to depend on~$p$. Specifically, we denote by~$a_p$ and~$b_p$ the  coefficients in the BEKK model when the dimension is~$p$. 

\subsection{Limiting property of ~ESD of sample covariance matrix under large~BEKK model}\label{LSD_BEKK}

\subsubsection{Reducible case}

Note that if $a_p=0$,  then the BEKK model  reduces to the \mbox{i.i.d.} case with~$\boldsymbol{\Sigma}_t\equiv \overline{\boldsymbol{\Sigma}}$. In general, if~$a_p$ is close to $0$, then the~BEKK model will be similar to the \hbox{i.i.d.}~case. 

Recall that for any two distributions $F_1$ and $F_2$,  the Levy distance between them is defined as
$$L(F_1,F_2):=\inf\{\varepsilon >0|F_1(x-\varepsilon )-\varepsilon \leq F_2(x)\leq F_1(x+\varepsilon )+\varepsilon {\mathrm  {\;for\;all\;}}x\in {\mathbb  {R}}\}.$$
It is well known that convergence in Levy distance implies convergence in distribution. 

{Define 
\begin{equation}\label{reducible_cond_para}
\eta(a, b, p)=\frac{a}{1-a-b}\min\Big(\sqrt{p(1-a-b)},1\Big).
\end{equation}}
The next theorem shows convergence in Levy distance between~$F^{\mathbf{S}_n}$ and~$F^{\mathbf{S}_n^0}$ {when $\eta(a_p, b_p, p)\to 0$.}
\begin{thm}\label{reducable}
Under model \eqref{BEKK} and Assumption \ref{asump1}, if \hbox{$\eta(a_p, b_p, p)\to 0$} as \mbox{$p\to \infty$},  then
\begin{equation}\label{thm1_1}
L(F^{\mathbf{S}_n},  F^{\mathbf{S}_n^0})=o_p(1).
\end{equation}
\end{thm}
Theorem \ref{reducable} implies that when $\eta(a_p, b_p, p)\to 0$ as \mbox{$p\to \infty$}, the~ESD of sample covariance matrix under the~BEKK model converges to the same M-P law as the \hbox{i.i.d.} case. We refer to the case when $\eta(a_p, b_p, p)\to 0$ as the \emph{reducible case}. Under the~reducible case, the population spectrum can be recovered by reversing the Mar$\check{\text{c}}$enko-Pastur law.

\subsubsection{Non-reducible case}
When the reducible condition does not hold, what will the~ESD of the sample covariance matrix be like? We have seen in~Figure~\ref{Fig1} that when~\mbox{$a=0.05$} and~$b=0.9$, which is a typical setting calibrated from empirical findings [\cite{ELW19}], the~ESD under the BEKK model appears to be more heavy-tailed than the~\mbox{i.i.d.} case. We refer to the case when~{$\eta(a_p, b_p, p)$ is bounded away from zero as the \emph{non-reducible case}. In practice, the two coefficients $a_p$ and $b_p$ learned from financial data appear to fit the non-reducible case. Therefore, the investigation of the non-reducible case is not only of theoretical interest but also  practically relevant. We measure the difference in the~ESD's between the~BEKK model and the \mbox{i.i.d.} case by the second moment of the~ESD's
:
$$M_2=M_2^p=\sum_{i=1}^p \widehat{\lambda}_i^2/p=\tr\big((\mathbf{S}_n)^2\big)/p\text{, and }M^0_2=M^{0,p}_2=\sum_{i=1}^p(\widehat{\lambda}_i^0)^2/p=\tr\big((\mathbf{S}^0_n)^2\big)/p.$$
For the \mbox{i.i.d.} case,~$E(M_2^{0,p})=yH_1^2+H_2+o(1)$, where $H_1=\lim_{p\to\infty}\tr(\overline{\boldsymbol{\Sigma}})/p$ and $H_2=\lim_{p\to \infty}\tr(\overline{\boldsymbol{\Sigma}}^2)/p$; see (4.14) of \cite{yin1986limiting}.
The next theorem gives the property of~$E(M^p_2)$  when~$\eta(a_p, b_p, p)$ is bounded away from zero. 
\begin{thm}\label{nonreducable}
Under model \eqref{BEKK} and Assumption \ref{asump1}, if \mbox{$\eta(a_p, b_p, p)>c$} for some constant $c>0$, then there exists $\delta>0$ such that for all $p$ large enough,
\begin{equation}\label{case1_thm2}
E(M^p_2)\geq E(M_2^{0,p})+\delta.
\end{equation}
\end{thm}
Theorem \ref{nonreducable} implies that, under the non-reducible case, the~ESD under the BEKK models is more heavy-tailed than the~\mbox{i.i.d.} case, a feature that is suggested by Figure \ref{Fig1}.
\subsection{Time-variation adjusted spectrum estimator}\label{ADJ_LSD}

Theorems \ref{reducable} and \ref{nonreducable} suggest that, unlike the \mbox{i.i.d.} case, the usual spectrum estimator based on the M-P law  does not always work under the~BEKK model. To recover the population spectrum under the BEKK model, a new estimator needs to be developed.

Recovering the population spectrum under large~BEKK models can be done by establishing Mar$\check{\text{c}}$enko-Pastur type equations for the sample covariance matrix. However, this is  challenging  due to the nonlinear dependency in returns.  
In the present paper, we will not pursue this direction. Instead, we provide an alternative solution using a time-variation adjustment approach.

\cite{ZL11} study a similar problem under elliptical models. They propose a self-normalization approach to remove the time variation in the covariance matrices.  Motivated by this idea, we aim to adjust the dynamic volatilities of the nonlinearly dependent data so that the adjusted data behave asymptotically  \mbox{i.i.d.}. Compared with the elliptical model considered in \cite{ZL11}, removing the time-varying dependency in~BEKK models is more challenging. 

Under the~BEKK model~\eqref{BEKK}, the time-variation in the covariance matrix is governed by the two coefficients,~$a_p$ and~$b_p$, and each~$R_{it}$ follows a univariate~GARCH model: 
\begin{equation}\label{uni_garch}
\sigma^2_{i,t+1}=(1-a_p-b_p)\overline{\sigma}_i^2+a_p R_{it}^2+b_p\sigma_{i,t}^2,
\end{equation}
where~$\sigma_{i,t}^2=(\boldsymbol{\Sigma}_{t})_{ii}$, and~$\overline{\sigma}_{i}^2=(\overline{\boldsymbol{\Sigma}})_{ii}$ for~$1\leq i\leq p$. 
As a result,~$a_p$ and~$b_p$ can be estimated  without  knowing the whole unconditional covariance matrix. Specifically, we  randomly select~one variable, say,~$i_0$,
fit a univariate~GARCH model to~$(R_{i_0t})$ and get QMLE~$\widehat{a}$ and~$\widehat{b}$:
$$(\widehat{a},\widehat{b},\widehat{\overline{\sigma}}_{i_0})=\mathop{\text{argmax}}_{(a,b,\overline{\sigma}_{i_0})\in \Omega} \frac{1}{n}\Bigg(\sum_{t=1}^n \Big(\frac{R_{i_0t}^2}{\sigma_{i_0t}^2}+\log(\sigma_{i_0t}^2)\Big)\Bigg),$$
where~$\Omega=\{(a,b,\overline{\sigma}), 0\leq a\leq 1, 0\leq b\leq 1, a+b\leq 1-\delta, \delta\leq\overline{\sigma}^2<C\}$ for some positive constants~$\delta$ and~$C$. The QMLE of the univariate~GARCH model is consistent with convergence rate  $\sqrt{1/n}$; see, Theorems 2.1 and 2.2 of \cite{francq2004maximum}. 
\begin{rem}In practice, we can use multiple variables to obtain a pooled estimator of $a_p$ and $b_p$. Numerical results suggest that a pooled estimator gives  a more stable estimation of $a_p$ and $b_p$. However, in terms of convergence rate, a pooled estimator does not necessarily help due to the cross-sectional dependence and the error in estimating the nuisance parameter $\overline{\sigma}_{i_0}$; see \cite{engle09} and \cite{pakel2011nuisance}. 

\end{rem}

We then use~$\widehat{a}$,~$\widehat{b}$ and past returns to construct a projection matrix:
\begin{equation}\label{P_mp}
\mathbf{P}_{t}=\frac{1-\widehat{a}-\widehat{b}+\widehat{a}\widehat{b}^{M_p}}{1-\widehat{b}}\mathbf{I}+\sum_{j=1}^{M_p} \widehat{a} \widehat{b}^{j-1} \mathbf{R}_{t-j}\mathbf{R}_{t-j}^T,
\end{equation}
where~$M_p$ represents the number of lagged returns included in $\mathbf{P}_t$, which grows with~$p$ and {$M_p=o(\sqrt{p})$.  The intuition behind such a definition is that, $\boldsymbol{\Sigma}_{t}=(1-a_p-b_p)/(1-b_p)\overline{\boldsymbol{\Sigma}}+\sum_{j=1}^\infty a_p b_p^{j-1} \mathbf{R}_{t-j}\mathbf{R}_{t-j}^T$, hence an appropriate choice of $M_p$ will make $\mathbf{P}_t^{-1/2}\boldsymbol{\Sigma}_t\mathbf{P}_t^{-1/2}$ close to $\overline{\boldsymbol{\Sigma}}$. 

To be more precise, using the projection matrix~$\mathbf{P}_t$, we define the time-variation adjusted  returns~$\widetilde{\mathbf{R}}_{t}=\mathbf{P}_{t}^{-1/2}\mathbf{R}_t$, and the time-variation adjusted (TV-adj) sample covariance matrix:
\begin{equation}\label{R_Mp}
\widetilde{\mathbf{S}}_n=\frac{1}{n}\sum_{t=1}^n\widetilde{\mathbf{R}}_{t}(\widetilde{\mathbf{R}}_{t})^T.
\end{equation}

\begin{thm}\label{recovered}
Under model \eqref{BEKK} and Assumption \ref{asump1}, 
if, in addition,~$\delta<\min(a_p, b_p)<a_p+b_p<1-\delta$ 
 for some~$\delta>0$, and $M_p$ satisfies that $M_p\to\infty$ and $M_p=o(\sqrt{p})$, then
\begin{equation}\label{thm3_1}
L(F^{\widetilde{\mathbf{S}}_n},  F^{\mathbf{S}_n^0})=o_p(1).
\end{equation}
\end{thm}

Theorem \ref{recovered} implies that the time-variation adjusted sample covariance matrix has the same~LSD as the \mbox{i.i.d.} case.  That is, under Assumption \ref{asump1}(iv), 
$F^{\widetilde{\mathbf{S}}_n}\overset{\text{P}}\to F$, and $F$ is determined by~$H$ in that its Stieltjes transform
\begin{equation}\label{Stieltjes}
m_F(z):=\int_{\lambda\in \mathbb{R}}\frac{1}{\lambda-z}dF(\lambda), z\in \mathbb{C}^+:=\{z\in \mathbb{C}, \im (z)>0\}
\end{equation}
solves the following equation
\begin{equation}\label{solve}
m_F(z)=\int_{\tau\in\mathbb{R}}\frac{1}{\tau\Big(1-y\big(1+zm_F(z)\big)\Big)-z}dH(\tau).
\end{equation}
See, e.g., Theorem~1 of \cite{MP67}.

We can then consistently estimate the population spectrum by reversing the Mar$\check{\text{c}}$enko-Pastur law.
Specifically, we denote the eigenvalues of~$\widetilde{\mathbf{S}}_n$ by~$\widetilde{\lambda}_1\geq ...\geq \widetilde{\lambda}_p$. We first regularize 
the eigenvalues of~$\widetilde{\mathbf{S}}_n$ to be~$\widetilde{\lambda}^\tau_i=\min (\widetilde{\lambda}_{i}, L)$ for some large constant~$L$. 
We then apply the Quantized Eigenvalues Sampling Transform (QuEST) algorithm in \cite{LW15} on~$(\widetilde{\lambda}^\tau_i)'s$ and obtain the estimated population spectrum. 
Denote by~$\widehat{\lambda}^{H}_1\geq \widehat{\lambda}^H_2\geq ...\geq \widehat{\lambda}^{H}_p$ the estimated eigenvalues and~$\lambda^H_1\geq \lambda^H_2\geq ...\geq \lambda^H_p$ the eigenvalues of~$\overline{\boldsymbol{\Sigma}}$.  
\begin{cor}\label{NLS_ADJ}Under the assumptions of~Theorem~\ref{recovered}, 
 $$ \frac{1}{p}\sum_{i=1}^p(\lambda^H_i-\widehat{\lambda}^H_i)^2=o_p(1).$$
\end{cor}
Corollary \ref{NLS_ADJ} guarantees that QuEST applied to the~TV-adj sample covariance matrix consistently estimates the population spectrum of the unconditional covariance matrix under large~BEKK models.

\subsection{Time-variation adjusted nonlinear shrinkage estimator of unconditional covariance matrix}\label{NLS_COV_TV}

The nonlinear shrinkage (NLS) estimator  [\cite{LW12,LW15,LW20}]  is structure-free and consistent in estimating the asymptotically optimal shrinkage estimator for \mbox{i.i.d.} data. In financial applications, the NLS has gained popularity in large portfolio optimization; see, e.g., \cite{LW17} and \cite{DLW21}. 

Motivated by the NLS developed under the~\mbox{i.i.d.} case, to estimate the unconditional covariance matrix under large BEKK models, we make use of the~TV-adj sample covariance matrix and consider rotation-equivariant  shrinkage estimators in the form~$\widehat{\boldsymbol{\Sigma}}=\sum_{i=1}^p\widehat{d}_i\widetilde{u}_i\widetilde{u}_i^T,$ 
where~$(\widetilde{u}_i)_{1\leq i\leq p}$ are eigenvectors of the TV-adj sample covariance matrix~$\widetilde{\mathbf{S}}_n$. The optimal rotation-equivariant estimator finds~$(\widehat{d}_i)_{1\leq i\leq p}$ that minimize  
$\|\widehat{\boldsymbol{\Sigma}}-\overline{\boldsymbol{\Sigma}}\|_{F}$.
 Elementary algebra shows that the optimal solution is $
\widehat{d}_i^*=\widetilde{u}_i^T\overline{\boldsymbol{\Sigma}}\widetilde{u}_i.$

In search of the asymptotically optimal shrinkage formula under large~BEKK models,  we study the following generalized empirical spectral distribution of the~\mbox{TV-adj} sample covariance matrix
\begin{equation}\label{F_general}
F^{\widetilde{\mathbf{S}}_n,g(\overline{\boldsymbol{\Sigma}})}(x)=\frac{1}{\tr\big(g(\overline{\boldsymbol{\Sigma}})\big)}\sum_{i=1}^p\Big(\widetilde{u}_i^{T}g(\overline{\boldsymbol{\Sigma}})\widetilde{u}_i\Big)\cdot\mathbf{1}_{[\widetilde{\lambda}_i,+\infty)}(x),
\end{equation}
which generalizes the~ESD of~$\widetilde{\mathbf{S}}_n$ by replacing the weight~$1/p$ with~$\widetilde{u}_i^Tg(\overline{\boldsymbol{\Sigma}})\widetilde{u}_i/\tr(g(\overline{\boldsymbol{\Sigma}}))$ for some bounded function~$g(\cdot)$, and~$g(\overline{\boldsymbol{\Sigma}})=\sum_{i=1}^p g(\lambda^H_i) v_iv_i^T$, where~$v_i$'s are the eigenvectors of~$\overline{\boldsymbol{\Sigma}}$. The limit of the generalized ESD of the sample covariance matrix under the~\mbox{i.i.d.} case is obtained in \cite{LP11} and is used to derive the asymptotically optimal shrinkage estimator. Parallel to the \mbox{i.i.d.} case, we study the limiting property of the generalized ESD of the TV-adj sample covariance matrix via the following  generalized Stieltjes transform:
\begin{equation}\label{theta_g}
\Theta^g_n(z)=
 \frac{1}{p}\tr\Big((\widetilde{\mathbf{S}}_n-z\mathbf{I})^{-1}g(\overline{\boldsymbol{\Sigma}})\Big).
\end{equation}

The following theorem gives the limit of~$\Theta^g_n(z)$. 
\begin{thm}\label{thetag}
Under model \eqref{BEKK} and Assumption \ref{asump1}, if, in addition,~$\delta<\min(a_p, b_p)<a_p+b_p<1-\delta$ for some~$\delta>0$, $M_p\to \infty$, $M_p=o(\sqrt{p})$,  the limiting distribution $H$ is supported by $[h_1, h_2]$ for some constants $0 < h_1\leq h_2 < \infty$, and $g$ is a bounded function on~$[h_1, h_2]$ with finitely many points of discontinuity,  then,$$\Theta^{g}_n(z)-\Theta^g(z)=o_p(1), \text{ for all }z\in \mathbb{C}^+,$$ where
\begin{equation}\label{limit_thetag}
 \Theta^g(z)=\int_{-\infty}^{+\infty}\Bigg(\tau\Big(1-y^{-1}-y^{-1} z m_F(z)\Big)-z\Bigg)^{-1} g(\tau) dH(\tau).
 \end{equation}
\end{thm}
The function $\Theta^g(z)$ is the limit of the generalized Stieltjes transform under the~\mbox{i.i.d.} case; see Theorem 2 of \cite{LP11}. Theorem \ref{thetag} states that the generalized~ESD based on the time-variation adjusted sample covariance matrix converges to the same limit as the \mbox{i.i.d.} case.
Therefore, we can utilize the same nonlinear shrinkage algorithm that is developed for~\mbox{i.i.d.} case to obtain the  time-variation adjusted nonlinear shrinkage estimator under~BEKK models.

Specifically, to estimate the unconditional covariance matrix, we perform the  nonlinear shrinkage algorithm by \cite{LW15} on~$\widetilde{\mathbf{S}}^{\tau}_n$, where~$\widetilde{\mathbf{S}}^{\tau}_n=\sum_{i=1}^p \widetilde{\lambda}_i^{\tau}\widetilde{u}_i\widetilde{u}_i^T$, $\widetilde{\lambda}_i^{\tau}=\min(\widetilde{\lambda}_i, L)$, and~$L$ is  a large constant.
{The truncation 
is applied to ensure that the support of the~ESD is bounded. 
We denote by~$\widetilde{\boldsymbol{\Sigma}}$ the resulting  covariance matrix estimator, which we call the~TV-adj nonlinear shrinkage estimator (TV-adj NLS).   Define~$\widetilde{\boldsymbol{\Sigma}}^{or}=\sum_{i=1}^pd_i^{or}(\widetilde{\lambda}^{\tau}_i)\widetilde{u}_i\widetilde{u}_i^T$, where
\begin{equation}\label{dioracle}
d_i^{or}(\widetilde{\lambda}^\tau_i)=
\begin{cases}
\frac{1}{(y-1) \check{m}_{\underline{F}}(0)} &\text{ if } \widetilde{\lambda}^\tau_i=0\text{ and }y>1,\\
\frac{\widetilde{\lambda}^\tau_i}{\big(1-y-y\widetilde{\lambda}^
\tau_i\cdot \check{m}_F(\widetilde{\lambda}^\tau_i)\big)^2}&\text{otherwise},
\end{cases}
\quad \text{ for }i=1, ..., p,
\end{equation}
$m_{\underline{F}}=(y-1)/z+ym_F(z)$,~$\underline{F}(x)=(1-y)\mathbf{1}_{\{[0,\infty)\}}(x)+yF(x)$,~$\check{m}(\lambda)=\lim_{z\in \mathbb{C}^+\to \lambda} m_F(z)$, $F(\cdot)$ and~$m_F(z)$ are given in \eqref{Stieltjes} and \eqref{solve}, respectively. By~Theorem~\ref{thetag} and Theorem~4 of \cite{LP11},~$\widetilde{\boldsymbol{\Sigma}}^{or}$ is the infeasible oracle shrinkage estimator.

\begin{thm}\label{NLS_ADJ_cov}Under the assumptions of~Theorem~\ref{thetag}, 
$$\frac{1}{\sqrt{p}}\|\widetilde{\boldsymbol{\Sigma}}-\widetilde{\boldsymbol{\Sigma}}^{or}\|_F=o_p(1).$$
\end{thm}
Theorem \ref{NLS_ADJ_cov} guarantees that the~\mbox{TV-adj NLS} consistently estimates the oracle shrinkage estimator under large~BEKK models in terms of convergence in dimension-normalized Frobenius norm. The convergence result achieved by the~\mbox{TV-adj NLS} under~BEKK models matches with that of the ordinary NLS under the \mbox{i.i.d.} case; see Proposition~4.3 of \cite{LW12} and  Theorem~3.1 of \cite{LW15}.

\section{Simulation Studies}\label{Simu}
\subsection{Simulation setup}
We generate data from the~BEKK model \eqref{BEKK} with~\mbox{$\mathbf{z}_t\underset{\text{i.i.d.}}\sim {N}(0, \mathbf{I})$}. The unconditional covariance matrix is set to be~$\overline{\boldsymbol{\Sigma}}=(\rho^{|i-j|})_{1\leq i,j\leq p}$, where\footnote{We have also evaluated other settings of unconditional covariance matrix including~$\rho=0, 0.2, 0.6$ and~$0.8$. 
The results are qualitatively similar.}~$\rho=0.4$. The dimension is set to be~\mbox{$p=100$} or $500$. We fix~$p/n=0.8$.  
About the parameters $(a,b)$:
\begin{itemize}
\item First, we choose four representative cases:
$(a, b)\in\{(0,0), (0.15,0.25),(0.1,0.65)$,
$(0.05,0.9)\}$. The setting~$(a,b)=(0,0)$ corresponds to the \mbox{i.i.d.} case, which is presented as a benchmark, and the other~$(a,b)$ pairs correspond to nontrivial~BEKK cases sorted 
 with increasing magnitudes of~$\eta(a, b, p)$ defined in \eqref{reducible_cond_para}, representing increasing levels of deviation from the \mbox{i.i.d.} case. The last configuration~$(a,b)=(0.05,0.9)$ is the setting used in \cite{ELW19} and is calibrated from empirical financial returns.  We simulate under each setting 100 replications and present the results in Section \ref{selected_ab}.

\item Next, we  examine more choices of $(a,b)$. Specifically,  we consider a grid of $(a,b)$'s in the region $\{(a,b):0.05\leq a\leq 0.5, 0.05\leq b\leq 0.90, a+b\leq 0.95\}$ and report the average results from 100 replications in Section \ref{ab_comprehensive}.
\end{itemize}


\subsection{Simulation results}
\subsubsection{Four $(a,b)$ cases}\label{selected_ab}
In this subsection, we present the simulation results for four representative $(a,b)$ cases, $(a, b)\in\{(0,0), (0.15,0.25),(0.1,0.65)$, $(0.05,0.9)\}$. 
\vskip 0.2cm
{\noindent\bf{Empirical spectral distribution of the sample covariance matrices}}
\vskip 0.1cm
We  compute the original sample covariance matrix~$\mathbf{S}_n=\sum_{t=1}^n{\mathbf{R}}_t{\mathbf{R}}_t^T/n$ and the~TV-adj sample covariance matrix~$\widetilde{\mathbf{S}}_n=\sum_{t=1}^n\widetilde{\mathbf{R}}_t\widetilde{\mathbf{R}}_t^T/n$, and compare their~ESDs with that of  the sample covariance matrix under the \mbox{i.i.d.} case,~$\mathbf{S}_n^0=\sum_{t=1}^n\mathbf{R}^0_{t}(\mathbf{R}^0_{t})^{T}/n$. 

We first illustrate the~ESDs from one random realization for~\mbox{$p=500$} in~Figure~\ref{Fig_BEKK_ESD}. We see that for all four cases, the~ESDs of the~TV-adj sample covariance matrix match remarkably well with that of the \mbox{i.i.d.} case. On the other hand,  under nontrivial~BEKK models,  the~ESDs of the original sample covariance matrices deviate from the M-P law, in particular, are more heavy-tailed.

\begin{figure}[H]
\begin{center}
\includegraphics[width=0.45\textwidth]{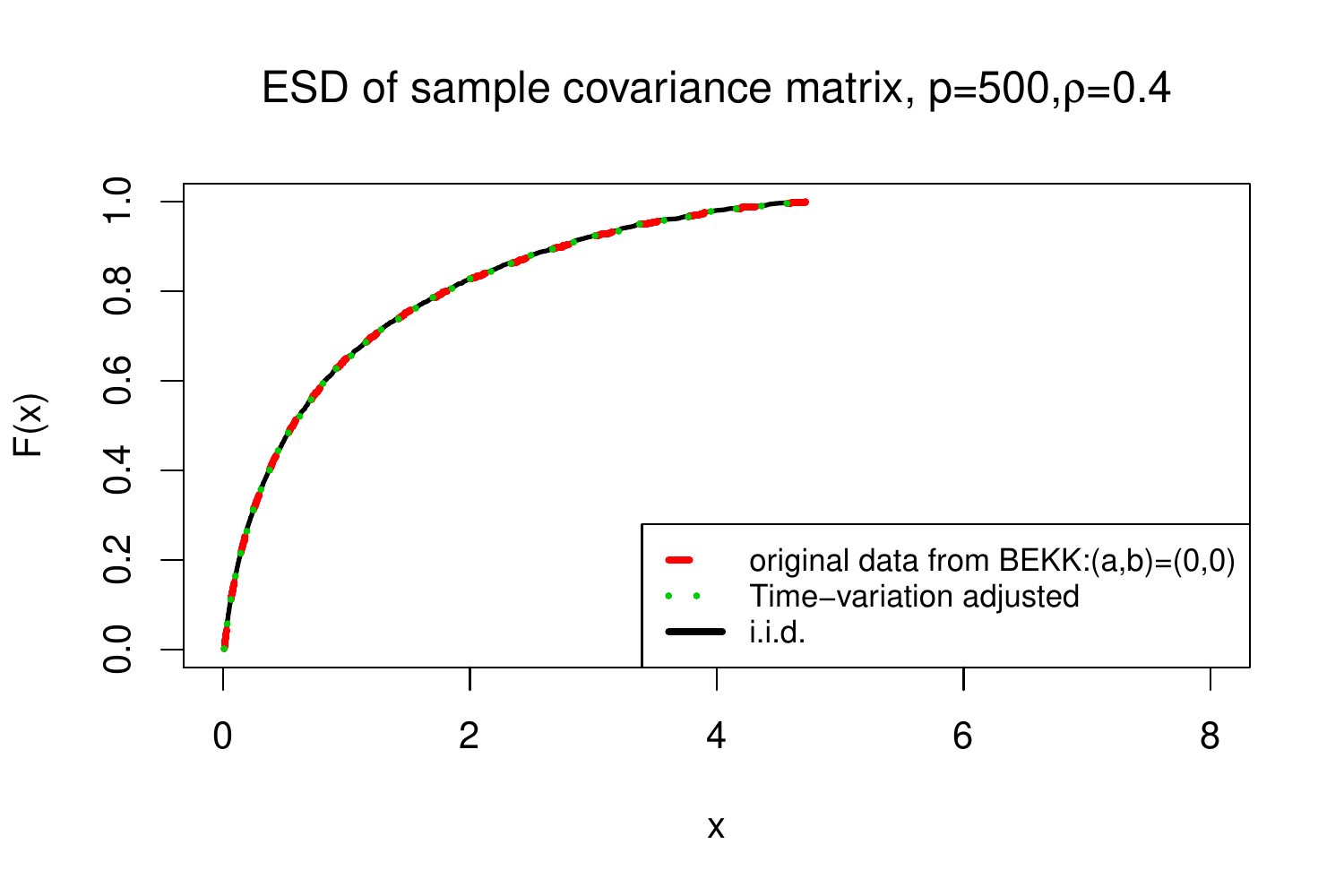}
\includegraphics[width=0.45\textwidth]{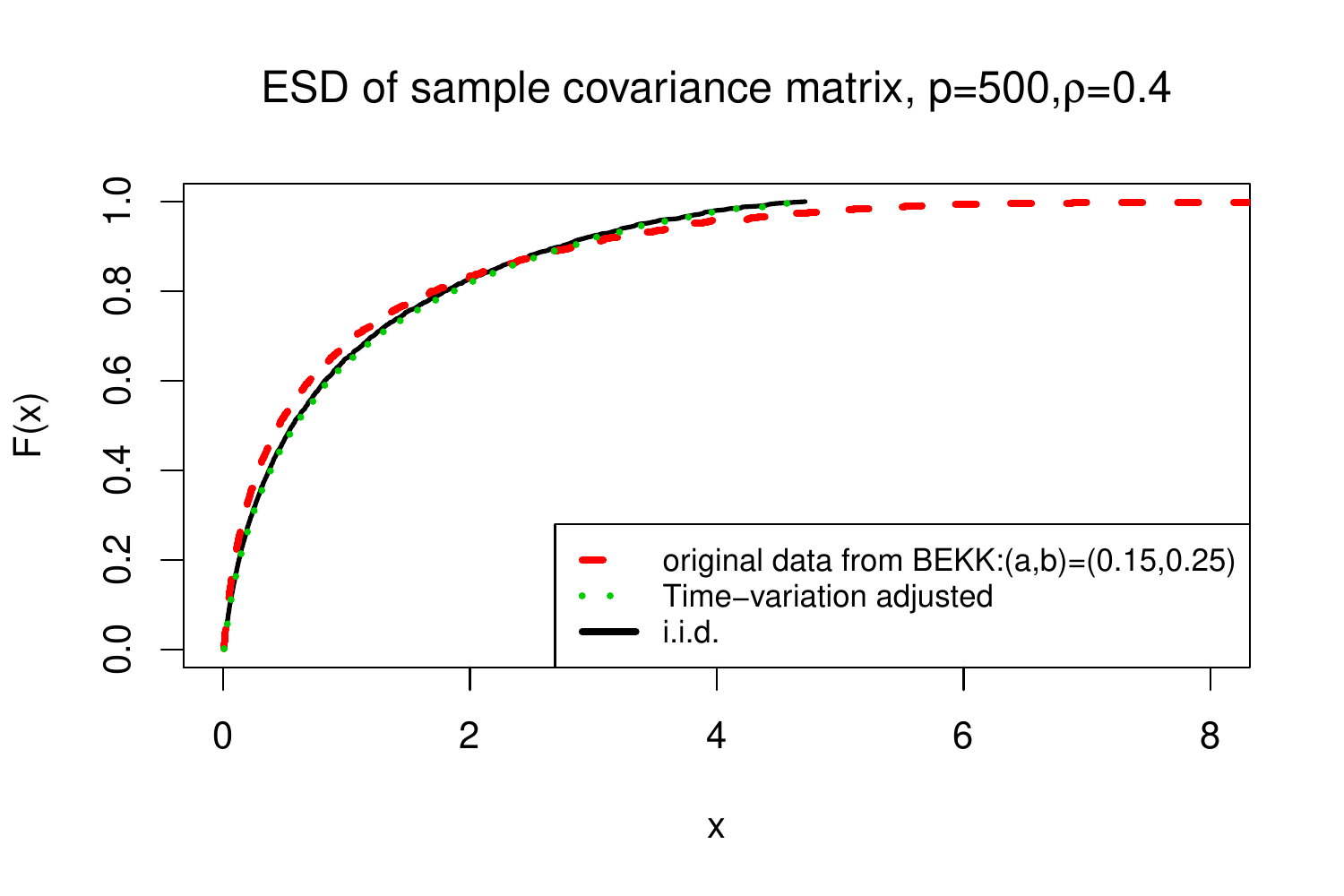}
\includegraphics[width=0.45\textwidth]{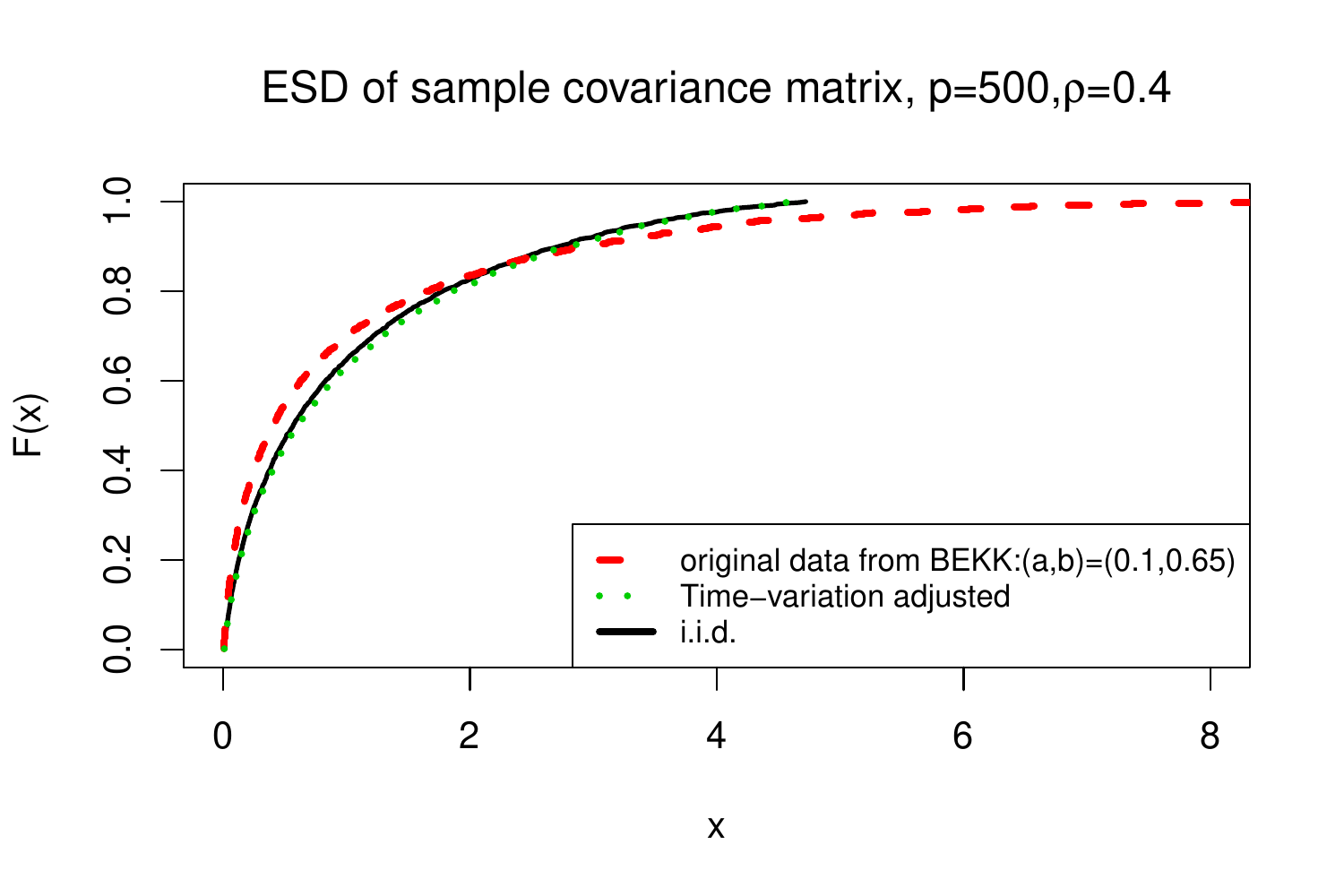}
\includegraphics[width=0.45\textwidth]{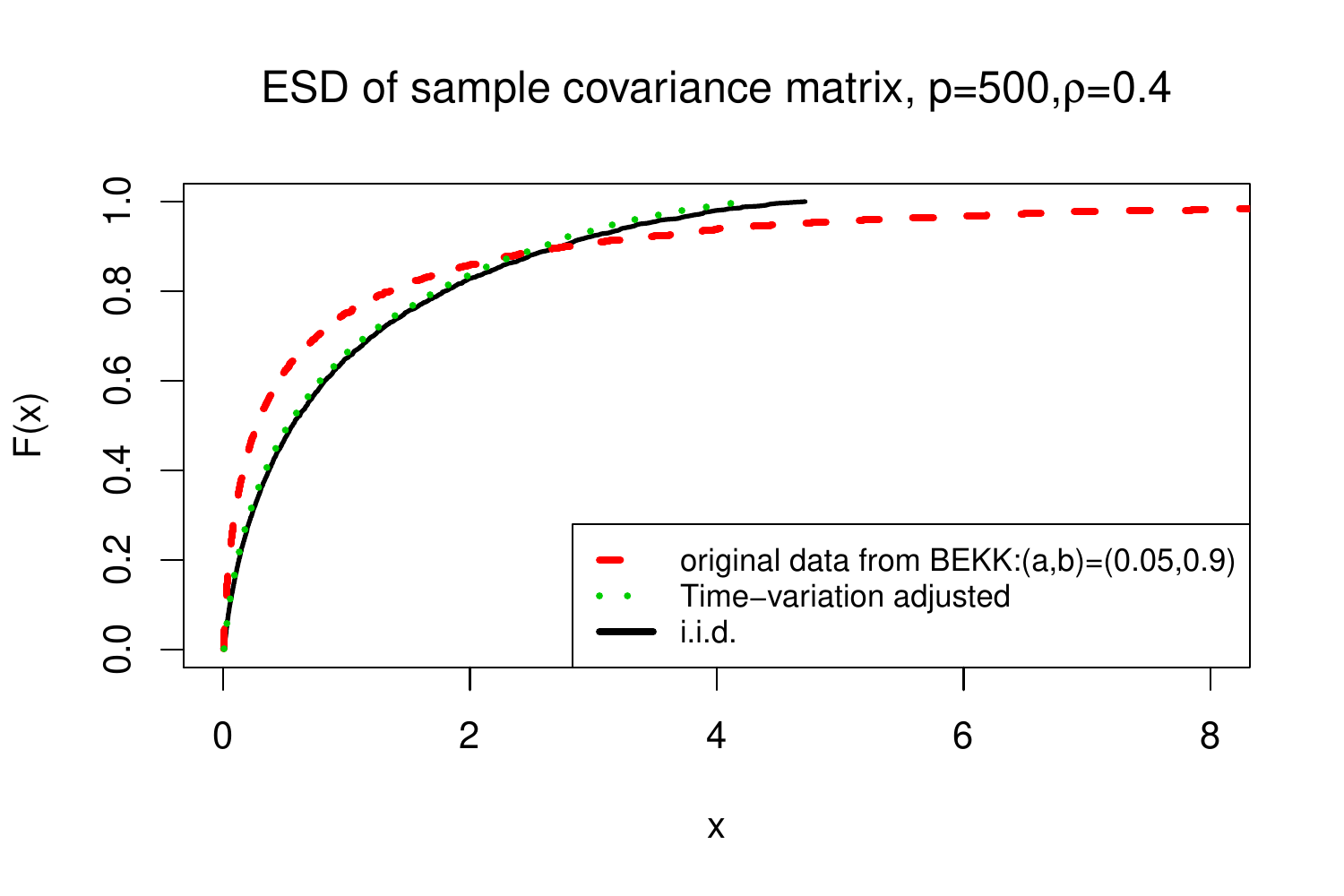}
\end{center}
\caption{{ESDs of sample covariance matrices of the original data,  the time-variation adjusted data, and the  \mbox{i.i.d.} data for~$p=500$,~$n=625$. The unconditional covariance matrix~$\overline{\boldsymbol{\Sigma}}=(0.4)^{|i-j|}$;~$(a,b)\in\{(0,0), (0.15,0.25), (0.1,0.65), (0.05,0.9)\}$.  
}}\label{Fig_BEKK_ESD}
\end{figure}

We then perform 100 replications and summarize the  Euclidean distance between the eigenvalues of~$\widetilde{\mathbf{S}}_n$ and~${\mathbf{S}}^0_n$, and between the eigenvalues of~$\mathbf{S}_n$ and~${\mathbf{S}}_n^0$ in~Table~\ref{simu_esd}. We see that:
\begin{itemize}
\item Under~nontrivial BEKK settings, the distance between the~ESD of the original sample covariance matrix and that under the \mbox{i.i.d.} case  increases  with increasing magnitudes of~$\eta(a, b, p)$. 
\item The distance between the ESD of the~TV-adj sample covariance matrix and that under the \mbox{i.i.d.} case is smaller and closer to zero under various nontrivial BEKK settings. Moreover, it decreases as $p$ gets larger.   
\item Comparing the performance of the TV-adj sample covariance matrix across different $(a,b)$ settings, the distance is the largest for the fourth setting $(a,b)=(0.05,0.9)$, in which case~$a+b$ is close to one, and the serial dependence is high. 
\end{itemize}

\begin{table}[H]
\caption{{{\it Summary of the distance~$\sqrt{\sum_{1\leq i\leq p}(\widehat{\lambda}_{i}-\widehat{\lambda}_{i}^0)^2}$, where~$(\widehat{\lambda}_{i}^0)_{1\leq i\leq p}$ are eigenvalues of~$\mathbf{S}^0_n$, and~$(\widehat{\lambda}_{i})_{1\leq i\leq p}$ are eigenvalues of~${\mathbf{S}}_n$ or~$\widetilde{\mathbf{S}}_n$. 
The mean and standard deviation (in parenthesis) from~100 replications are reported.
}}}
\begin{center}\small
\tabcolsep 0.08in\renewcommand{\arraystretch}{1} \doublerulesep
2pt
\begin{tabular}{clccc}
\hline \hline
$(a,b)$&&(0.15,0.25)&(0.1,0.65)&(0.05,0.9) \\ \hline
\multirow{3}{*}{$(p,n)=(100,125)$}&&\\
&$\mathbf{S}_n$&0.277&0.413&0.907\\
&&(0.078)&(0.106)&(0.195)\\
&$\widetilde{\mathbf{S}}_n$&0.089&0.123&0.215\\
&&(0.034)&(0.033)&(0.079)\\
&&&\\
\multirow{3}{*}{$(p,n)=(500,625)$}&&\\
&$\mathbf{S}_n$&0.279&0.429&1.046\\
&&(0.045)&(0.050)&(0.120)\\
&$\widetilde{\mathbf{S}}_n$&0.029&0.054&0.162\\
&&(0.017)&(0.029)&(0.056)\\
\hline

\end{tabular}
\end{center}\label{simu_esd}
\end{table}

\vskip 0.2cm
{\noindent\bf{Population spectrum estimation}}

\vskip 0.1cm
Next, we evaluate the estimators of the population eigenvalues. The performance of the proposed time-variation adjusted NLS spectrum estimator\footnote{The function ``tau\_estimate'' from R package ``nlshrink'' is used to compute the estimated eigenvalues. } (\mbox{TV-adj NLS}-Spectrum)  is compared with that of the NLS spectrum estimator based on the original sample covariance matrix (original NLS-Spectrum). 
We measure the estimation error by 
$\sqrt{\sum_{1\leq i\leq p}(\lambda_i^{H}-\widehat{\lambda}_i^{H})^2}$, where~$\widehat{\lambda}_1^{H}\geq \widehat{\lambda}_2^{H}\geq ...\geq \widehat{\lambda}_p^{H}$ are the estimated eigenvalues, and~$\lambda_1^{H}\geq \lambda_2^{H}\geq ...\geq \lambda_p^{H}$ are the eigenvalues of~$\overline{\boldsymbol{\Sigma}}$. 

 In~Figure~\ref{Fig_BEKK_LSDH}, we plot the distributions of the estimated eigenvalues from one random realization with~$p=500$. 
We see that, under all four cases, the proposed~TV-adj spectrum estimator is close to the population spectrum, and its performance is similar to  the spectrum estimator based on the infeasible \mbox{i.i.d.} data. On the other hand, the shrinkage spectrum estimator based on the original sample covariance matrix significantly deviates  from the population spectrum.
\begin{figure}
\begin{center}
\includegraphics[width=0.45\textwidth]{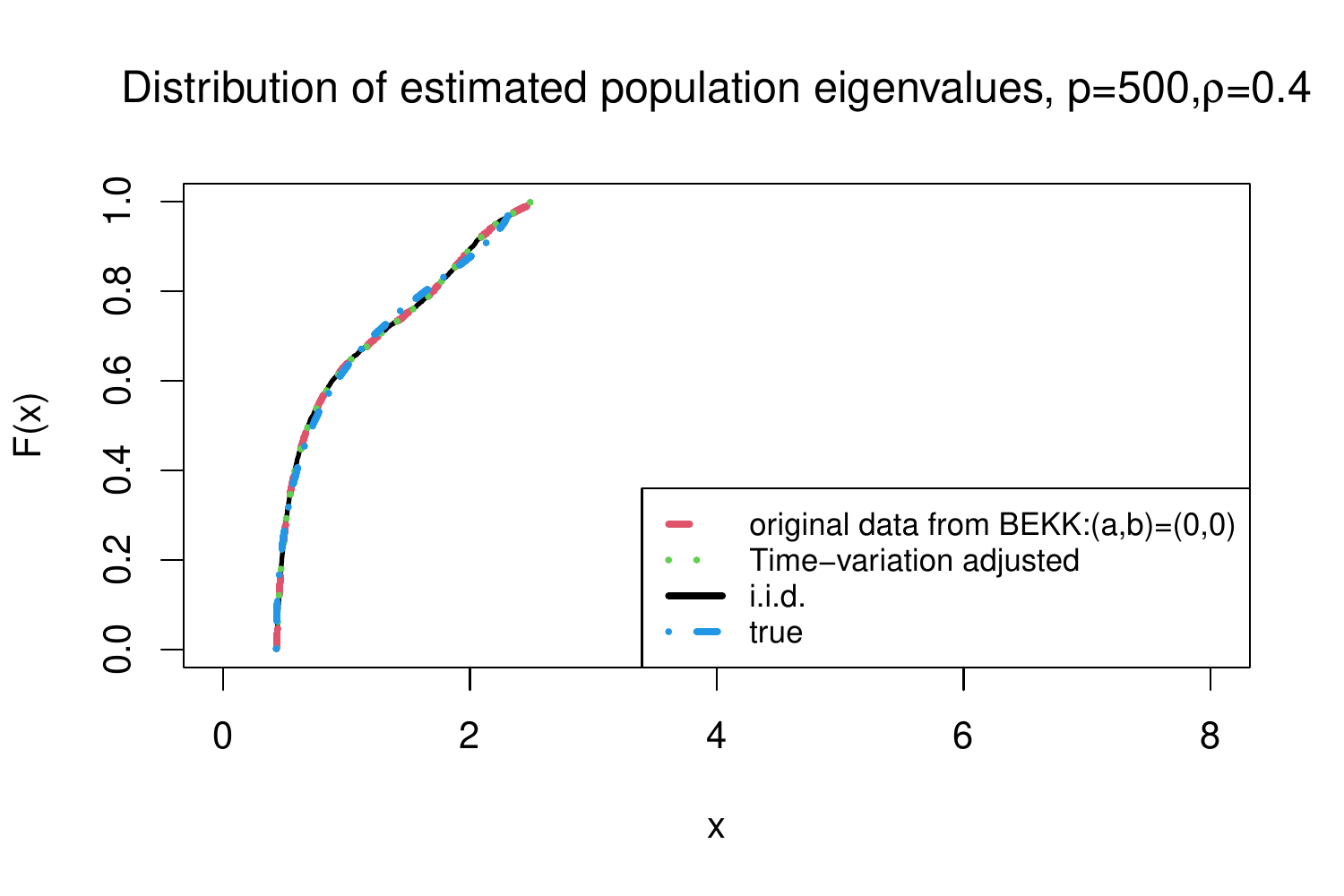}
\includegraphics[width=0.45\textwidth]{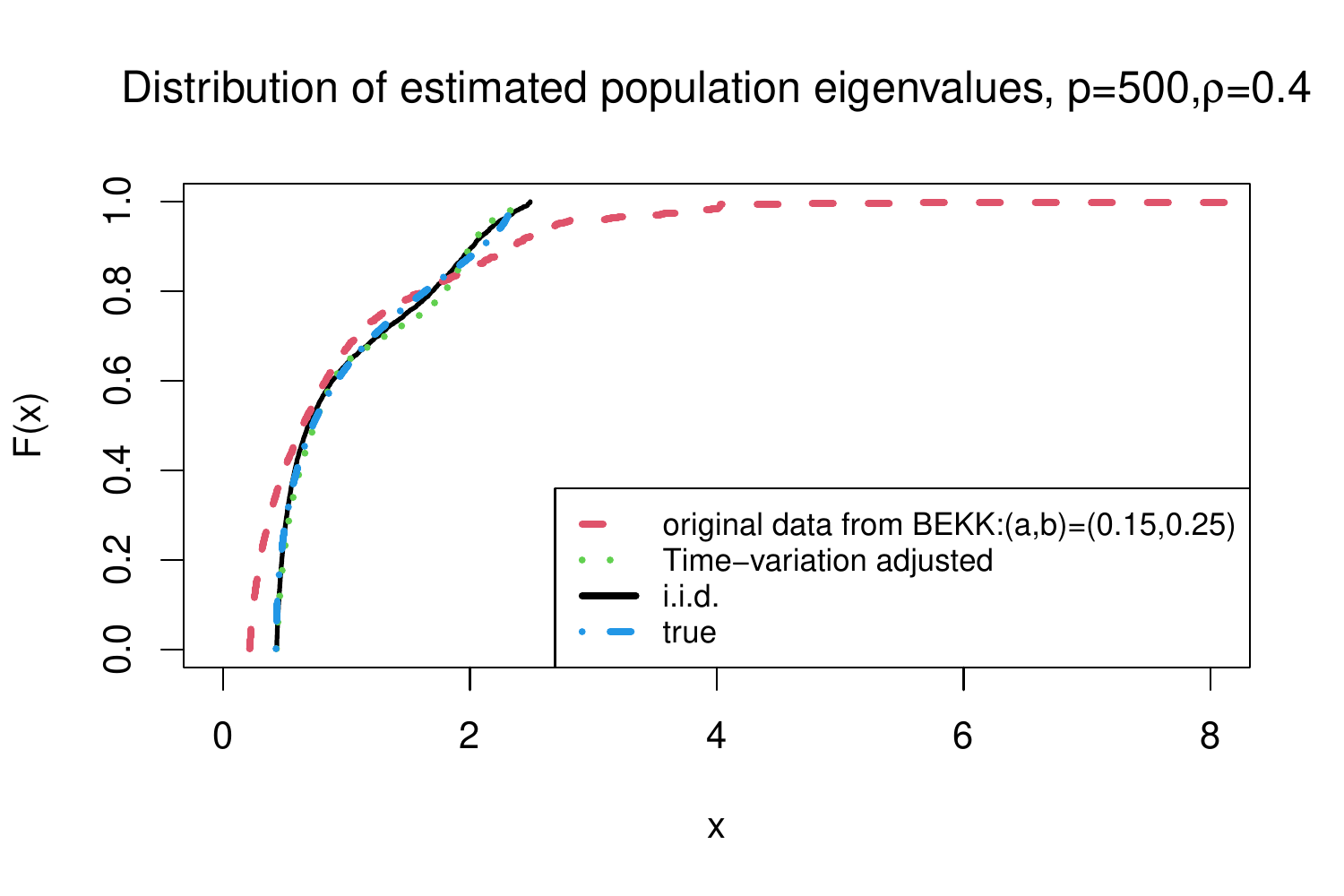}
\includegraphics[width=0.45\textwidth]{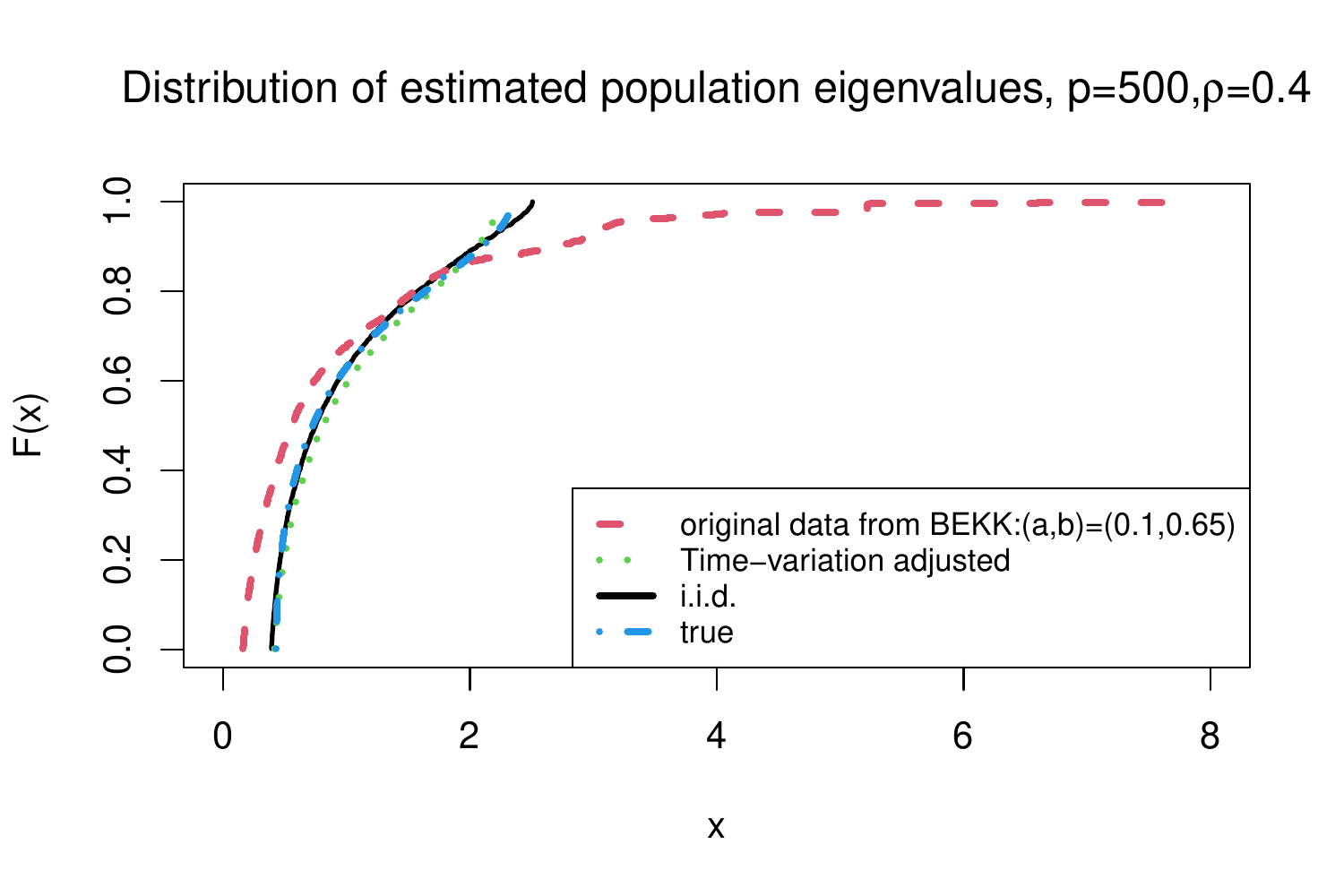}
\includegraphics[width=0.45\textwidth]{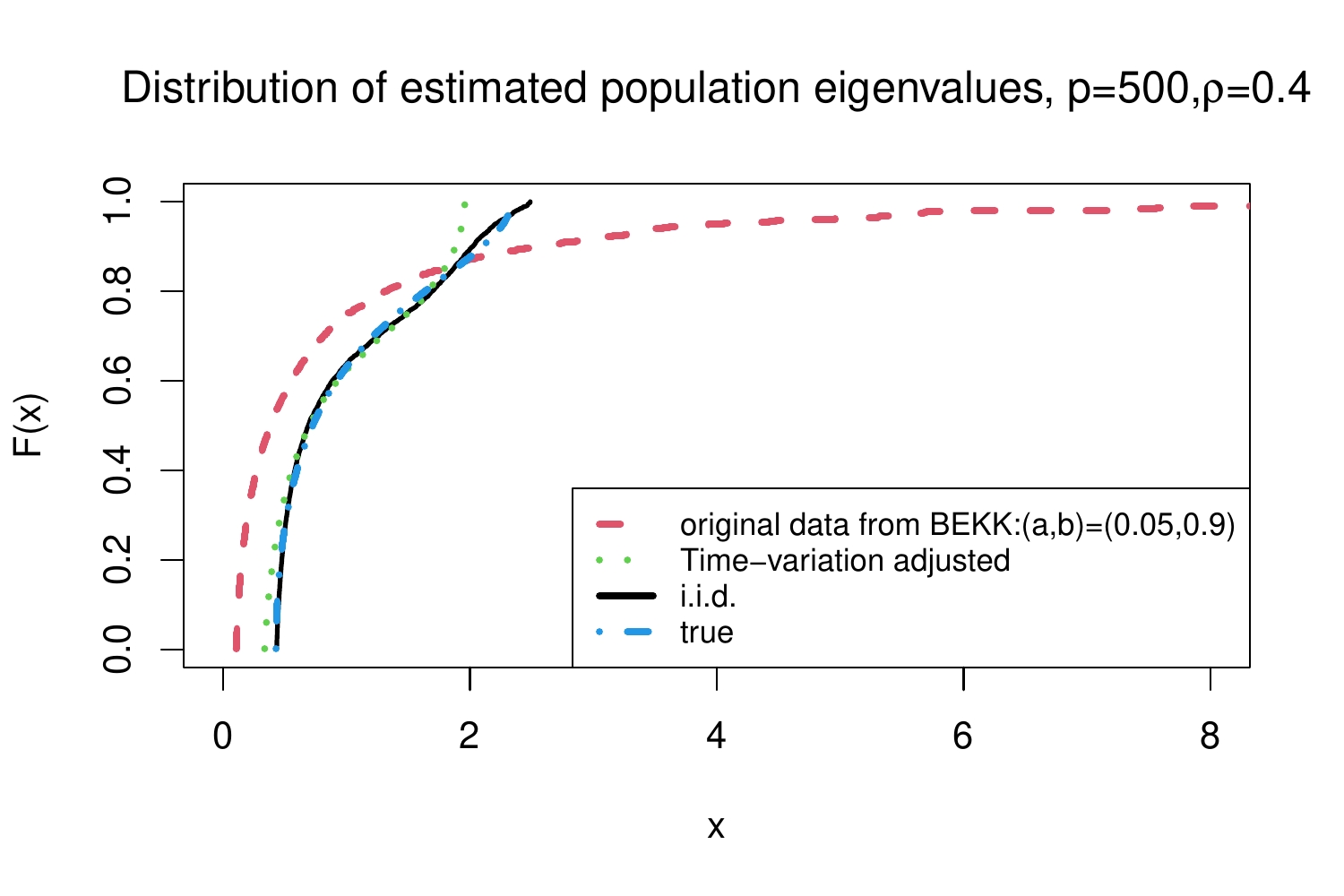}
\end{center}
\caption{{Distributions of estimated eigenvalues. They are obtained by nonlinear shrinkage estimators applied to the original data  the time-variation adjusted data, and the \mbox{i.i.d.} data. The dimension~$p=500$,~$n=625$.~$\overline{\boldsymbol{\Sigma}}=(0.4)^{|i-j|}$; $(a,b)\in\{(0,0), (0.15,0.25), (0.1,0.65), (0.05,0.9)\}$.  
}}\label{Fig_BEKK_LSDH}
\end{figure}

The Euclidean distances between the estimated population eigenvalues and the true ones from~100 replications are summarized in~Table~\ref{simu_LSDH}. We see that:

\begin{itemize}
\item The error of the original NLS-spectrum estimator increases with~$\eta(a,b,p)$. It also gets larger as the dimension gets higher.

\item The proposed TV-adj NLS-spectrum estimator  dominantly outperforms the original NLS-spectrum estimator with a substantially lower estimation error. 

\item As the dimension $p$ grows, the performance of the TV-adj NLS-spectrum estimator gets closer to the performance of the shrinkage estimator under the \mbox{i.i.d.} case.

\item The performance of the TV-adj NLS-spectrum estimator is only slightly worse than the infeasible shrinkage estimator based on~\mbox{i.i.d.} data for the first two nontrivial BEKK settings.  
For the fourth setting when $(a,b)=(0.05,0.9)$, because $a+b$ is close to one, the estimation error of the TV-adj NLS-spectrum estimator is larger. However, as the dimension $p$ grows, the error decreases  and becomes closer to that of the shrinkage estimator under the~\mbox{i.i.d.} case. 
\end{itemize}

\begin{table}[H]
\caption{{{\it Summary of the distance between the estimated population eigenvalues and the true ones. The mean and standard deviation (in parenthesis) from~100 replications
 are reported. The pair~$(a,b)=(0,0)$ represents the \mbox{i.i.d.} case and is presented as the benchmark. The remaining~$(a, b)$ pairs are for nontrivial~BEKK cases.
}}}
\begin{center}\small
\tabcolsep 0.08in\renewcommand{\arraystretch}{1} \doublerulesep
2pt
\begin{tabular}{llc|ccc}
\hline \hline
&&\mbox{i.i.d.} &\multicolumn{3}{c}{BEKK}\\
&($a$,$b$)&(0,0) &(0.15,0.25)&(0.1,0.65)&(0.05,0.9) \\ \hline
\multirow{4}{*}{$(p,n)=(100,125)$}&&\\
&original NLS-Spectrum&0.136&0.411&0.566&1.109\\
&&(0.041)&(0.094)&(0.124)&(0.209)\\
&TV-adj NLS-Spectrum&0.143&0.217&0.251&0.313
\\
&&(0.047)&(0.071)&(0.071)&(0.100)\\
&&&\\
\multirow{4}{*}{$(p,n)=(500,625)$}&&\\
&original NLS-Spectrum &0.052&0.401&0.577&1.243\\
&&(0.026)&(0.049)&(0.053)&(0.119)\\
&TV-adj NLS-Spectrum&0.054&0.057&0.088&0.184
\\
&&(0.027)&(0.025)&(0.026)&(0.035)\\
\hline
\end{tabular}
\end{center}\label{simu_LSDH}
\end{table}
\vskip 0.2cm
{\noindent\bf{Unconditional covariance matrix estimation}}

\vskip 0.1cm

Finally, we evaluate the unconditional covariance matrix estimation. We compute the NLS estimators\footnote{The function ``nlshrink\_cov'' in R package ``nlshrink'' is used in computing the nonlinear shrinkage estimator of the covariance matrix.} based on the time-variation adjusted sample covariance matrix (\mbox{TV-adj NLS}) and the original sample covariance matrix (original NLS). 
The estimation error is measured by the Frobenius norm~$\sqrt{\sum_{1\leq i,j\leq p}( \overline{\boldsymbol{\Sigma}}_{ij}-\widehat{\boldsymbol{\Sigma}}_{ij})^2}$, where~$\widehat{\boldsymbol{\Sigma}}$ is the estimated unconditional covariance matrix. The results are summarized in Table~\ref{simu_Sigma}.
\begin{table}[H]
\caption{{{\it Summary statistics of the estimation error of the estimated unconditional covariance matrix in Frobenius norm~$\sqrt{\sum_{1\leq i,j\leq p}(\widehat{\boldsymbol{\Sigma}}_{ij}-\overline{\boldsymbol{\Sigma}}_{ij})^2}$. The mean and standard deviation (in parenthesis) from~100 replications are reported. The pair~$(a,b)=(0,0)$ corresponds to the \mbox{i.i.d.} case and is presented as the benchmark. The remaining~$(a, b)$ pairs correspond to nontrivial~BEKK cases.
}}}
\begin{center}\small
\tabcolsep 0.08in\renewcommand{\arraystretch}{1} \doublerulesep
2pt
\begin{tabular}{llc|ccc}
\hline \hline
&&\mbox{i.i.d.} &\multicolumn{3}{c}{BEKK}\\
&($a$,$b$)&(0,0) &(0.15,0.25)&(0.1,0.65)&(0.05,0.9) \\ \hline
\multirow{4}{*}{$(p,n)=(100,125)$}&&\\
&original NLS&5.079&6.668&7.933&12.810\\
&&(0.048)&(0.558)&(0.894)&(1.824)\\
&TV-adj NLS&5.080&5.219&5.308&6.433\\
&&(0.049)&(0.096)&(0.078)&(0.711)\\
&&\\
\multirow{4}{*}{$(p,n)=(500,625)$}&&\\
&original NLS&11.314&14.850&17.878&31.384\\
&&(0.021)&(0.664)&(0.866)&(2.405)\\
&TV-adj NLS&11.320&11.362&11.469&11.970\\
&&(0.026)&(0.028)&(0.050)&(0.196)\\
\hline
\end{tabular}
\end{center}\label{simu_Sigma}
\end{table}
We see from Table \ref{simu_Sigma} that:
\begin{itemize}
\item The estimation error of the original NLS increases sharply as $\eta(a, b, p)$ gets large and as the dimension grows. 
\item The proposed~\mbox{TV-adj NLS} greatly improves over the original NLS with a substantially lower estimation error. 

\item The performance of the \mbox{TV-adj NLS} is only slightly worse than that of the NLS under the~\mbox{i.i.d.} case for the first two nontrivial $(a,b)$ settings.  For the most challenging case $(a,b)=(0.05,0.9)$, because $a+b$ is close to one, the error of \mbox{TV-adj NLS} is larger. However,  when $p$ grows, the performance becomes closer.
\end{itemize}

\subsubsection{Performance under more choices of $(a,b)$}\label{ab_comprehensive}

In this subsection, we present the simulation results for a grid of $(a,b)$'s in the region $\{(a,b):0.05\leq a\leq 0.5, 0.05\leq b\leq 0.90, a+b\leq 0.95\}$. 
\vskip 0.2cm
\noindent{\bf Empirical spectral distribution of the sample covariance matrices}
\vskip 0.1cm

In Figure \ref{SESD_distance_sphere}, we plot the average Euclidean distance between the eigenvalues of~$\mathbf{S}_n$ and~${\mathbf{S}}_n^0$, and between the eigenvalues of~$\widetilde{\mathbf{S}}_n$ and~${\mathbf{S}}^0_n$. We see that the Euclidean distance between the~ESD of the original sample covariance matrix and that under the \mbox{i.i.d.} case grows substantially as $a$ and $a+b$ increase. In contrast, the distance of the eigenvalues of the~TV-adj sample covariance matrix to that under the \mbox{i.i.d.} case is close to zero for various $(a,b)$ settings. The distance surface for the TV-adj sample covariance matrix is almost flat,  except when $a+b$ approaches one, but when the dimension $p$ increases,  it again becomes flatter and closer to zero. 

\begin{figure}[H]
\begin{center}
\includegraphics[width=0.48\textwidth]{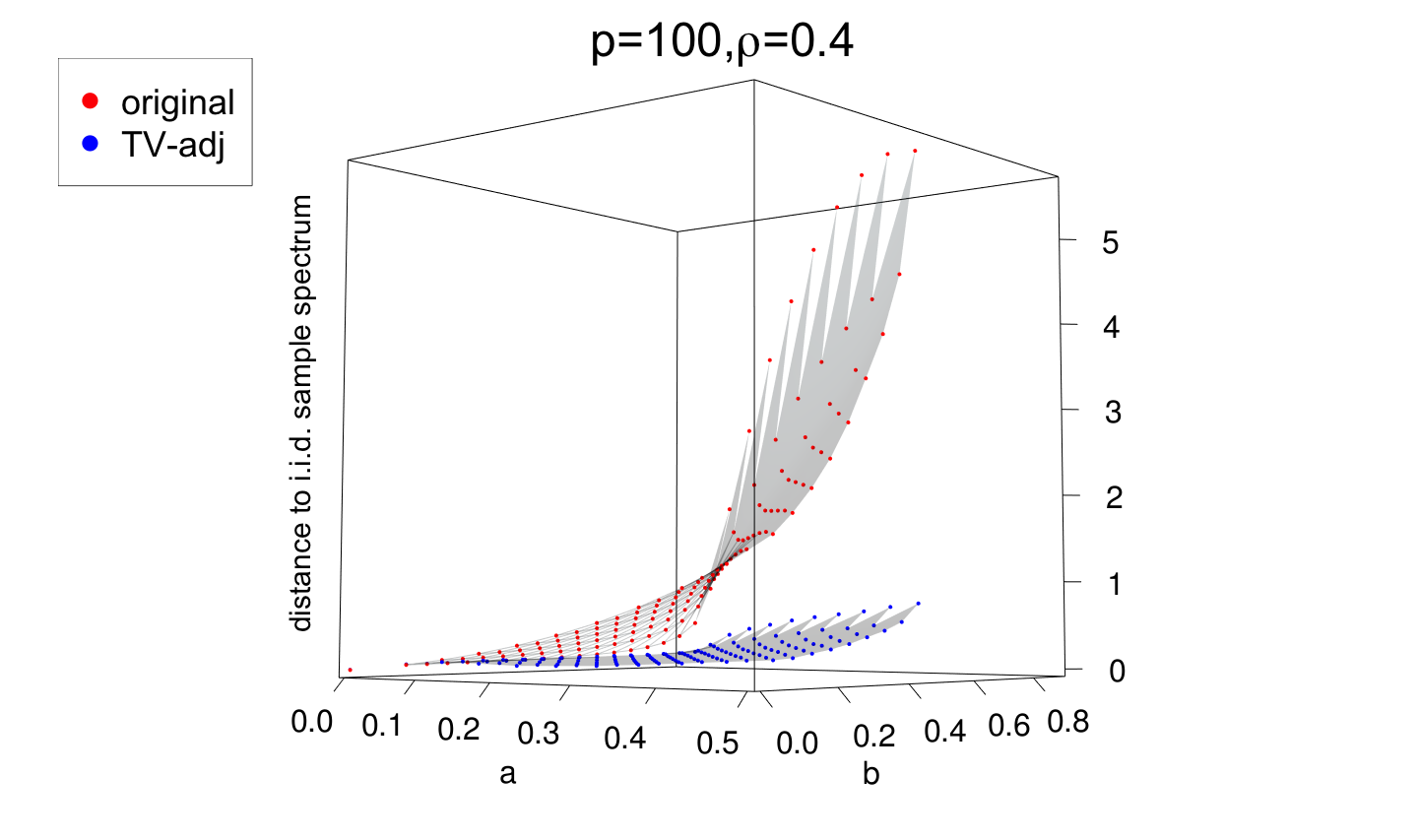}
\includegraphics[width=0.48\textwidth]{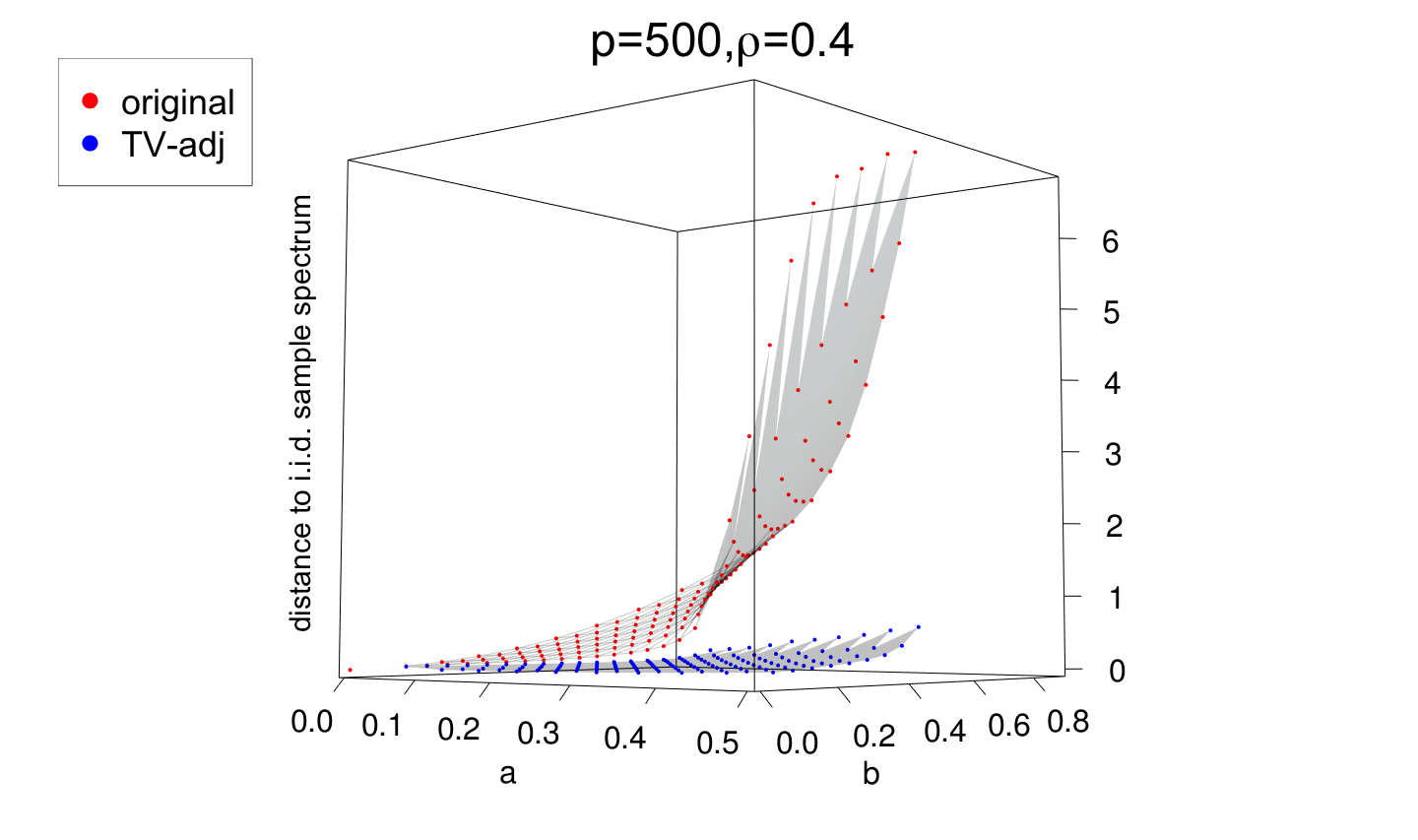}
\end{center}
\caption{{Euclidean distance between the eigenvalues of the sample covariance matrix/TV-adj sample covariance matrix under the BEKK model and the eigenvalues of the sample covariance matrix under the \mbox{i.i.d.} case for \mbox{$p=100$ (left)} and $p=500$ (right). The unconditional covariance matrix is~$\overline{\boldsymbol{\Sigma}}=(0.4)^{|i-j|}$. The evaluation is made for a grid of $(a,b)$'s in the region $\{(a,b):0.05\leq a\leq 0.5, 0.05\leq b\leq 0.90, a+b\leq 0.95\}$.
}}\label{SESD_distance_sphere}

\end{figure}
\vskip 0.2cm

\noindent{\bf Population spectrum estimation}

\vskip 0.1cm
In Figure \ref{PESD_distance_sphere}, we plot the average Euclidean distance between the estimated  eigenvalues and the true eigenvalues. The methods under comparison are the original NLS estimator and the proposed \mbox{TV-adj NLS} estimator. We see that the error of the original NLS estimator increases sharply as $a$ and $a+b$ increase. It also gets larger when the dimension is higher. In contrast, the \mbox{TV-adj NLS} performs robustly well for various $(a,b)$ settings and it dominantly outperforms the original NLS in all cases. The error surface for the \mbox{TV-adj NLS} is almost flat except when $a+b$ approaches one, but it gets closer to zero when $p$ grows. 

\begin{figure}[H]
\begin{center}
\includegraphics[width=0.48\textwidth]{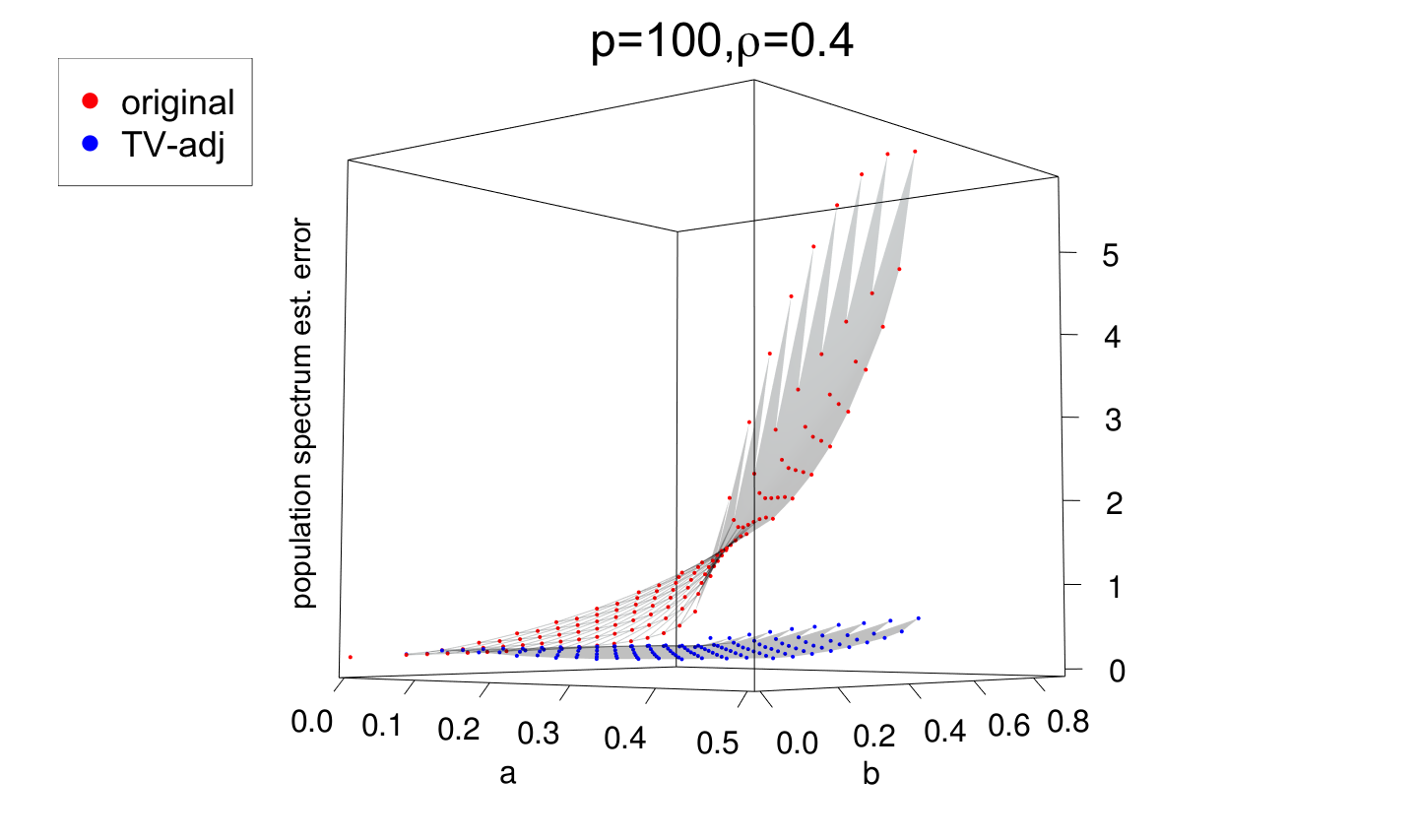}
\includegraphics[width=0.48\textwidth]{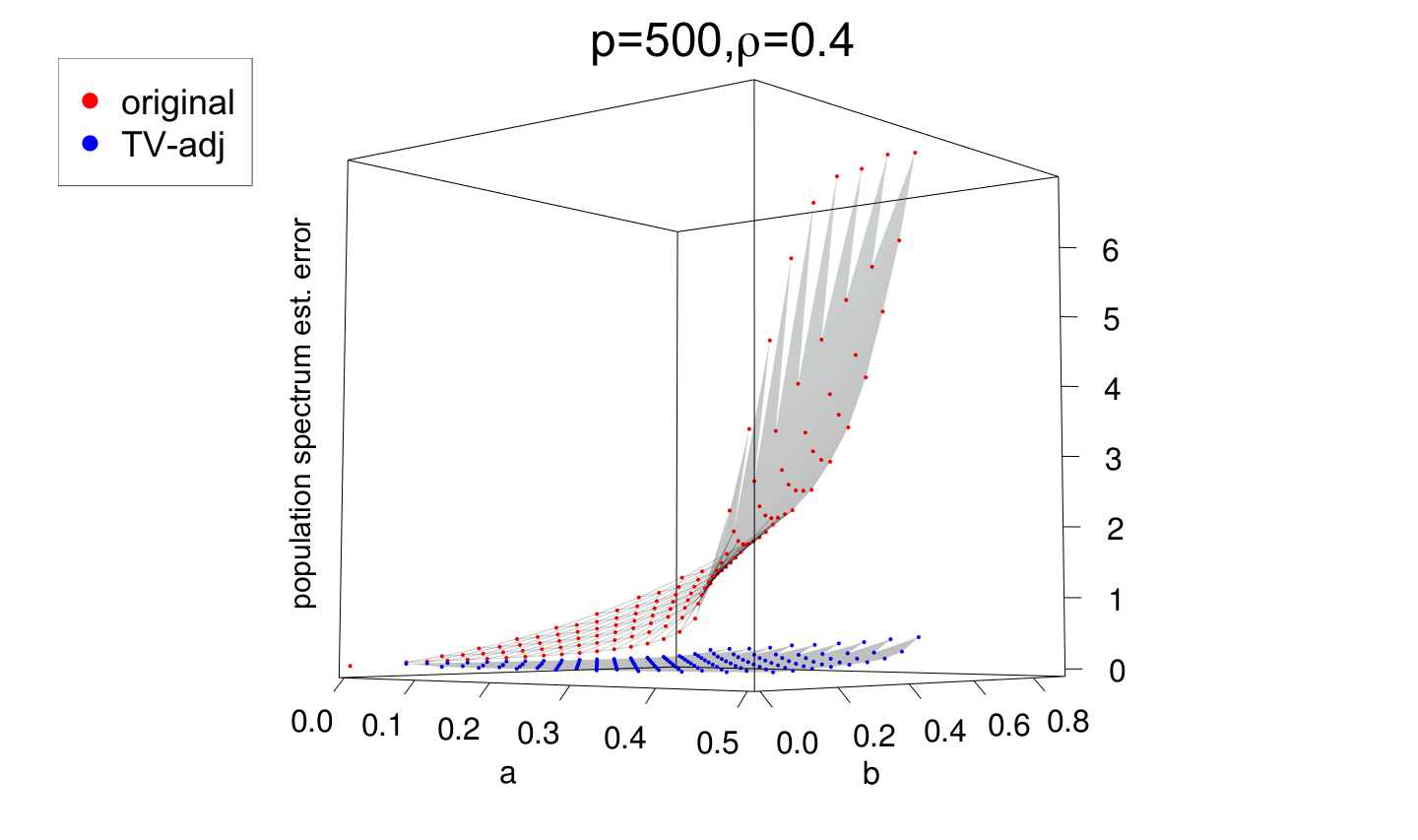}
\end{center}
\caption{{Estimation error of the population eigenvalues for \mbox{$p=100$ (left)} and $p=500$ (right). The unconditional covariance matrix is~$\overline{\boldsymbol{\Sigma}}=(0.4)^{|i-j|}$. The methods under comparison are the original NLS-spectrum estimator and the \mbox{TV-adj NLS}-spectrum estimator. The evaluation is made for a grid of $(a,b)$'s in the region $\{(a,b):0.05\leq a\leq 0.5, 0.05\leq b\leq 0.90, a+b\leq 0.95\}$.
}}\label{PESD_distance_sphere}

\end{figure}
\vskip 0.2cm

\noindent{\bf Unconditional covariance matrix estimation}
\vskip 0.1cm

Finally, in Figure \ref{CovEst_sphere}, we plot the average Frobenius error of the NLS and \mbox{TV-adj NLS} in estimating the unconditional covariance matrix. We see that the original NLS estimator performs poorly when $(a,b)$ deviates from $(0,0)$. When $p$ grows, the error also becomes bigger. The \mbox{TV-adj NLS} dominantly outperforms the NLS with a lower estimation error in all cases.  The error surface for the \mbox{TV-adj NLS} is almost flat and only  slightly higher near the edge when $a+b$ is close to one.

\begin{figure}[H]
\begin{center}
\includegraphics[width=0.48\textwidth]{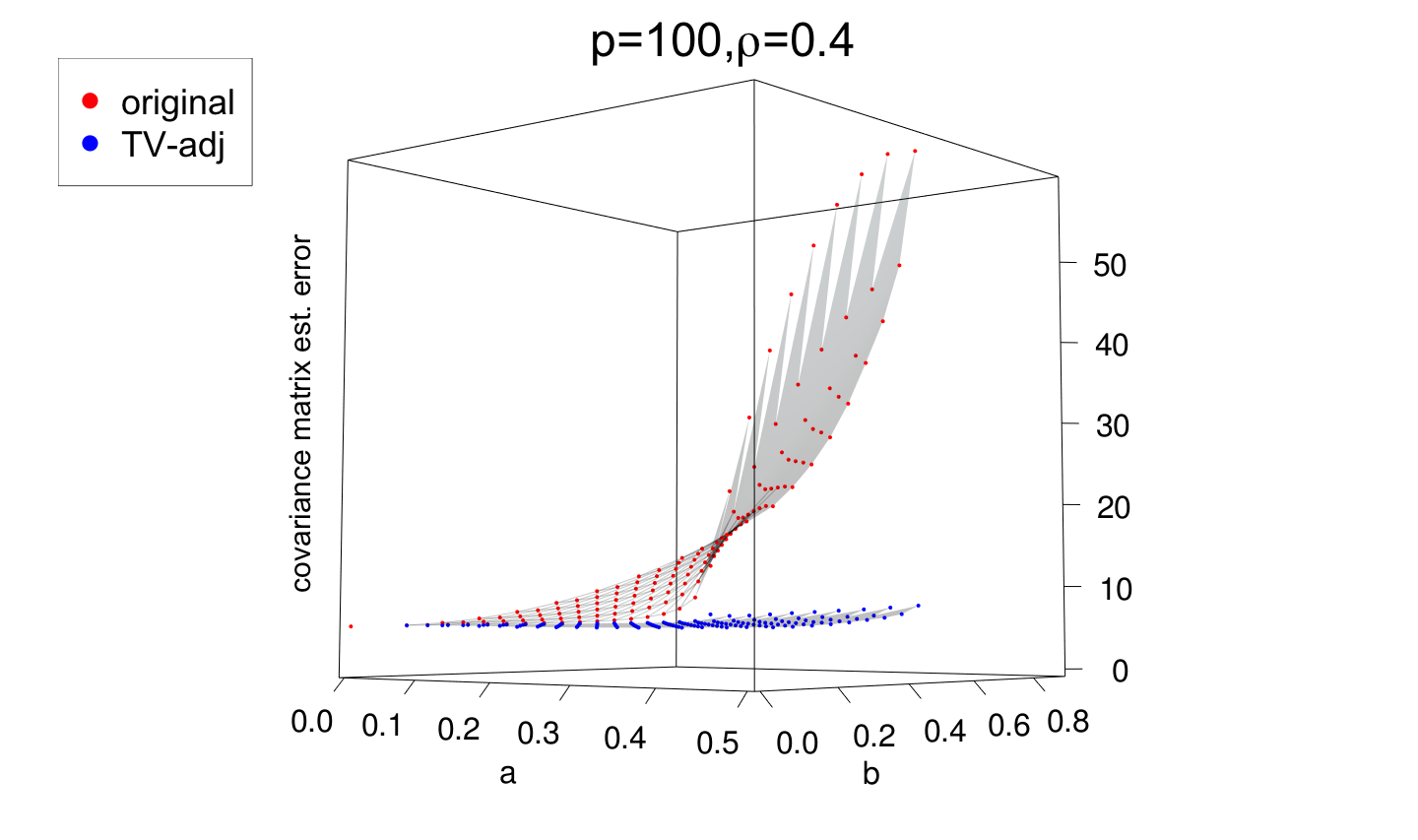}
\includegraphics[width=0.48\textwidth]{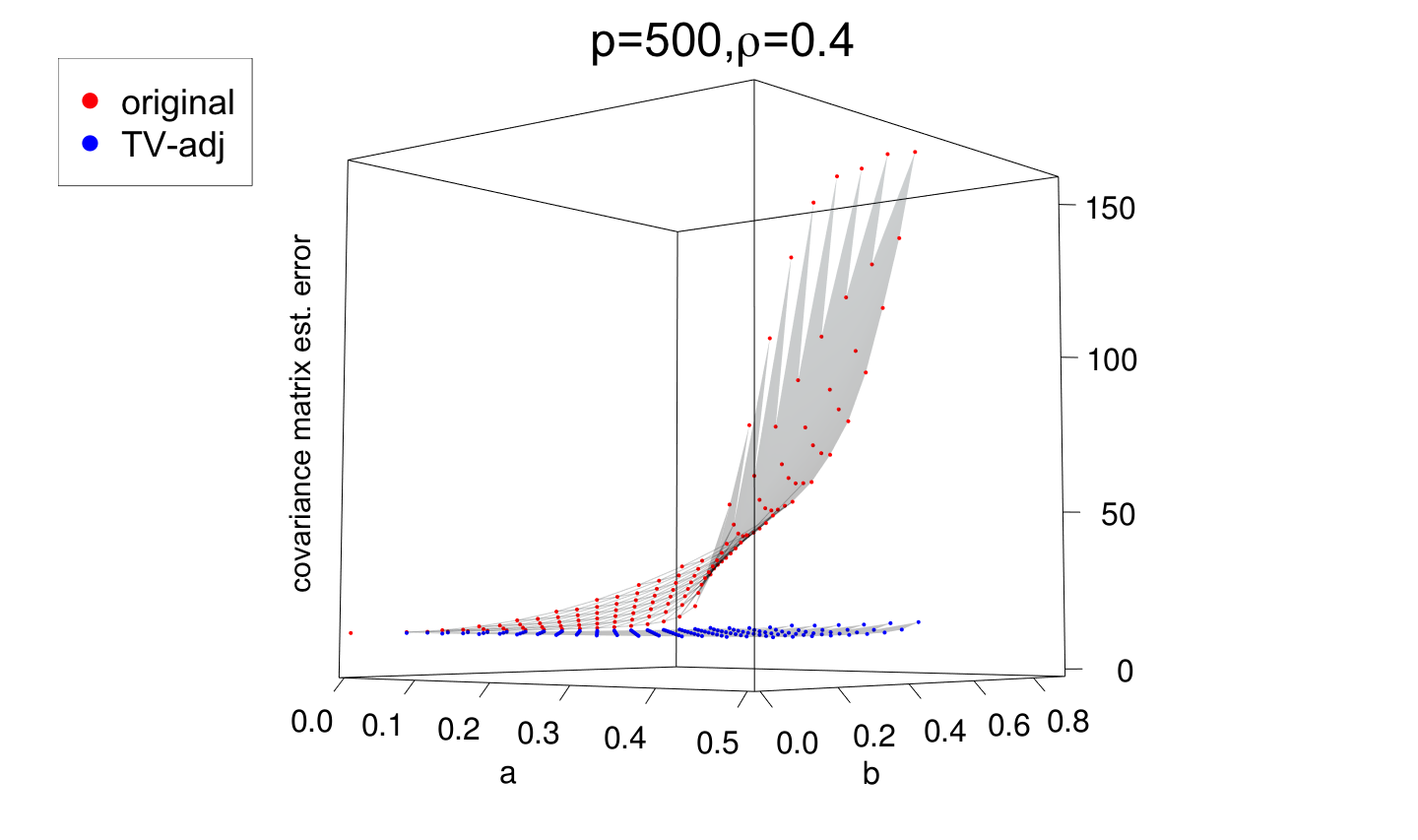}
\end{center}
\caption{{Estimation error of the unconditional covariance matrix in Frobenius norm. The unconditional covariance matrix is~$\overline{\boldsymbol{\Sigma}}=(0.4)^{|i-j|}$. The methods under comparison are the original NLS estimator and the \mbox{TV-adj NLS} estimator. The evaluation is made for a grid of $(a,b)$'s in the region $\{(a,b):0.05\leq a\leq 0.5, 0.05\leq b\leq 0.90, a+b\leq 0.95\}$.
}}\label{CovEst_sphere}

\end{figure}

\section{Conclusion}\label{Conc}

We investigate 
the limiting spectral properties of high-dimensional sample covariance matrices under dynamic volatility models. We show that under large~BEKK models, the asymptotics of the spectral distribution of sample covariance matrix depend on the asymptotic size of the~innovation coefficient and the persistence coefficient. In particular, we give explicit conditions under which the ESD has the same limit as or deviates from the \mbox{i.i.d.} case. 
 Furthermore, we develop a consistent estimator of the spectrum of the unconditional covariance matrix under large~BEKK models.  The proposed approach is based on a time-variation adjusted sample covariance matrix, for which we show that the LSD shares the same Mar$\check{\text{c}}$enko-Pastur law as the \mbox{i.i.d.} case. Finally, we 
 propose a~TV-adj nonlinear shrinkage estimator of the unconditional covariance matrix. The estimator is consistent in estimating the oracle shrinkage estimator under large~BEKK models. 
\section{Proof of Theorem \ref{recovered}}\label{proof_recovered}
In this section we prove Theorem \ref{recovered}.  The proofs of other main and technical results are given in the Supplementary Material \cite{DZ22_supp}.

We divide the proof of Theorem \ref{recovered} into three steps. In the first step, we show that replacing $(\widehat{a},\widehat{b})$ with $(a_p,b_p)$ in \eqref{P_mp} does not change the LSD of the TV-adj sample covariance matrices. In the second step, we show that the problem can be reduced into proving convergence in Frobenius norm for  the expected difference in the covariance matrices' square root matrics after some orthogonal transformation. In the last step, we show \eqref{E_tildR}.  For ease of notation, we drop the subscript~$p$ in~$a_p$ and~$b_p$ in the following proofs.

\underline{\emph{Step One}}:  Define 
$$
\aligned
&\check{\mathbf{P}}_{t}=\frac{1-a-b+ab^{M_p}}{1-b}\mathbf{I}+\sum_{j=1}^{M_p} a b^{j-1} \mathbf{R}_{t-j}\mathbf{R}_{t-j}^T,\\
&\check{{\mathbf{R}}}_t=(\check{\mathbf{P}}_{t})^{-1/2}\mathbf{R}_t,\text{ and}\\
&\check{{\mathbf{S}}}_n=\frac{1}{n}\sum_{t=1}^n\check{{\mathbf{R}}}_t\check{{\mathbf{R}}}_t^T.
\endaligned
$$
By Corollary A.42 of \cite{BS10}, 
\begin{equation}\label{L_check_diff}
L^4(F^{\widetilde{\mathbf{S}}_n},F^{\check{{\mathbf{S}}}_n})\leq \frac{2}{p}\tr (\widetilde{\mathbf{S}}_n+\check{{\mathbf{S}}}_n)\cdot \frac{1}{pn}\tr\Big((\widetilde{\mathbf{R}}-\check{{\mathbf{R}}})(\widetilde{\mathbf{R}}-\check{{\mathbf{R}}})^T\Big),
\end{equation}
where $\widetilde{\mathbf{R}}=(\widetilde{\mathbf{R}}_1, ...,\widetilde{\mathbf{R}}_n)$, and $\check{{\mathbf{R}}}=(\check{{\mathbf{R}}}_1, ...,\check{{\mathbf{R}}}_n)$. By definition,
$(\widehat{a}, \widehat{b})$ satisfies
\begin{equation}\label{QMLE_bound}
(\widehat{a}, \widehat{b},\widehat{\overline{\sigma}}_{i_0})\in \Omega=\{(a,b, \overline{\sigma}), 0\leq a\leq 1, 0\leq b\leq 1, a+b\leq 1-\delta, \delta\leq\overline{\sigma}<C\}.
\end{equation}
By the definition of $\mathbf{P}_t$ in \eqref{P_mp}, we have $\mathbf{P}_t\geq (1-\widehat{a}-\widehat{b})\mathbf{I}\geq \delta \mathbf{I}$. Similarly, under the assumption that $a+b<1-\delta$, we have $\check{\mathbf{P}}_t\geq (1-a-b)\mathbf{I}\geq\delta \mathbf{I}$. 
Define $$\widetilde{\boldsymbol{\Sigma}}_t=\mathbf{P}_t^{-1/2}\boldsymbol{\Sigma}_t\mathbf{P}_t^{-1/2}\text{, and }\,\check{{\boldsymbol{\Sigma}}}_t=\check{\mathbf{P}}_t^{-1/2}\boldsymbol{\Sigma}_t\check{\mathbf{P}}_t^{-1/2}.$$ Using the fact that if $\mathbf{A}\geq 0$ and $\mathbf{B}\geq 0$, then $\tr(\mathbf{A}\mathbf{B})\geq 0$,  
 we have that, for all $t$, 
\begin{equation}\label{bound_sigma_ii}
\tr({\widetilde{\boldsymbol{\Sigma}}}_t)\leq \frac{1}{\delta}\tr(\mathbf{\boldsymbol{\Sigma}}_t),\quad
\tr(\check{{\boldsymbol{\Sigma}}}_t)\leq \frac{1}{\delta}\tr(\mathbf{\boldsymbol{\Sigma}}_t),
\end{equation}
and
$$
\aligned
&\tr(\check{{\mathbf{R}}}_t\check{{\mathbf{R}}}_t^T)\leq \frac{1}{\delta}\tr(\mathbf{R}_t\mathbf{R}_t^T),\quad \tr(\widetilde{\mathbf{R}}_t\widetilde{\mathbf{R}}_t^T)\leq \frac{1}{\delta}\tr(\mathbf{R}_t\mathbf{R}_t^T). 
\endaligned$$
Therefore,
$$
\aligned
&\tr(\widetilde{\mathbf{S}}_n/p)\leq \frac{1}{\delta}\cdot\tr (\mathbf{S}_n/p),\quad \text{and}\quad \tr(\check{{\mathbf{S}}}_n/p)\leq \frac{1}{\delta}\cdot\tr (\mathbf{S}_n/p).
\endaligned
$$
By the independence between  $(\boldsymbol{\Sigma}_t)$ and $(\mathbf{z}_t)$ and Assumption \ref{asump1}(iii),  we have $E\Big(\tr (\mathbf{S}_n/p)\Big)=\tr(\overline{\boldsymbol{\Sigma}})/p=O(1)$. It follows that 
\begin{equation}\label{trS_ntilde}
\tr (\widetilde{\mathbf{S}}_n/p)=O_p(1),\quad \text{and }\; \tr (\check{{\mathbf{S}}}_n/p)=O_p(1).
\end{equation}
Write $$
\aligned\varepsilon_n=\max\Bigg(&\max_{1\leq j\leq M_p}\Big(\Big|\frac{ab^{j-1}}{\widehat{a}\widehat{b}^{j-1}}-1\Big|,\Big|\frac{\widehat{a}\widehat{b}^{j-1}}{ab^{j-1}}-1\Big|\Big),\\
&\Bigg|\frac{(1-\widehat{a}-\widehat{b}+\widehat{a}\widehat{b}^{M_p})(1-b)}{(1-{a}-{b}+{a}{b}^{M_p})(1-\widehat{b})}-1\Bigg|,\Bigg|\frac{(1-{a}-{b}+{a}{b}^{M_p})(1-\widehat{b})}{(1-\widehat{a}-\widehat{b}+\widehat{a}\widehat{b}^{M_p})(1-b)}-1\Bigg|\Bigg).
\endaligned$$
By Theorem 2.2 of \cite{francq2004maximum}, under the assumption that $\delta\leq a, b\leq a+b\leq 1-\delta$, we have \begin{equation}\label{con_rate}
\widehat{a}-a=O_{p}(1/\sqrt{n}),\;\text{ and  }\;\;\widehat{b}-b=O_{p}(1/\sqrt{n}).
\end{equation}
By the assumptions that $M_p=o(\sqrt{p})$, $\delta\leq a, b\leq a+b\leq 1-\delta$, $p/n\asymp 1$, and  \eqref{con_rate}, we get
\begin{equation}\label{e_n_rate}
\varepsilon_n=o_p(1).
\end{equation}
Note that when $\varepsilon_n<1$, we have
$$0<(1-\varepsilon_n)\check{\mathbf{P}}_t\leq\mathbf{P}_t\leq (1+\varepsilon_n)\check{\mathbf{P}}_t.$$
By the L$\ddot{\text{o}}$wner-Heinz inequality, 
$$\frac{1}{\sqrt{1+\varepsilon_n}}\check{\mathbf{P}}^{-1/2}_t\leq\mathbf{P}^{-1/2}_t\leq \frac{1}{\sqrt{1-\varepsilon_n}}\check{\mathbf{P}}^{-1/2}_t.
$$
By Weyl's theorem, we get that
\begin{equation}\label{error_both}
\Big\|\mathbf{P}^{-1/2}_t-\frac{1}{\sqrt{1+\varepsilon_n}}\check{\mathbf{P}}^{-1/2}_t\Big\|\leq \Big(\frac{1}{\sqrt{1-\varepsilon_n}}-\frac{1}{\sqrt{1+\varepsilon_n}}\Big) \cdot\|\check{\mathbf{P}}^{-1/2}_t\|. 
\end{equation}
Recall that $\check{\mathbf{P}}_t\geq \delta \mathbf{I}$, hence \begin{equation}\label{l2_boundp2}
\|\check{\mathbf{P}}^{-1/2}_t\|\leq \sqrt{\frac{1}{\delta}}.
\end{equation}
By
the triangle inequality, \eqref{e_n_rate}, \eqref{error_both} and \eqref{l2_boundp2},  we get that 
$$\aligned
\|\mathbf{P}^{-1/2}_t-\check{\mathbf{P}}^{-1/2}_t\|\leq&\Big\|\mathbf{P}^{-1/2}_t-\frac{1}{\sqrt{1+\varepsilon_n}}\check{\mathbf{P}}^{-1/2}_t\Big\|+ \Big(1-\frac{1}{\sqrt{1+\varepsilon_n}}\Big)\cdot\|\check{\mathbf{P}}^{-1/2}_t\|\\
=&O\Big(\frac{1}{\sqrt{1-\varepsilon_n}}-\frac{1}{\sqrt{1+\varepsilon_n}}+1-\frac{1}{\sqrt{1+\varepsilon_n}}\Big)\\
=&o_p(1).
\endaligned
$$
Moreover,
$$
\aligned
&\tr\Big((\widetilde{\mathbf{R}}_t-\check{{\mathbf{R}}}_t)(\widetilde{\mathbf{R}}_t-\check{{\mathbf{R}}}_t)^T\Big)\\
=&\mathbf{R}_t^T(\mathbf{P}_t^{-1/2}-\check{\mathbf{P}}_t^{-1/2})^2\mathbf{R}_t\\
\leq &\|\mathbf{P}_t^{-1/2}-\check{\mathbf{P}}_t^{-1/2}\|^2\cdot \|\mathbf{R}_t\|^2=\|\mathbf{P}_t^{-1/2}-\check{\mathbf{P}}_t^{-1/2}\|^2 \cdot(\mathbf{z}_t^T\boldsymbol{\Sigma}_t\mathbf{z}_t).
\endaligned
$$
Because $E(\mathbf{z}_t^T\boldsymbol{\Sigma}_t\mathbf{z}_t)=E\big(\tr(\boldsymbol{\Sigma}_t)\big)=\tr(\overline{\boldsymbol{\Sigma}})=O(p)$,  we have $\mathbf{z}_t^T\boldsymbol{\Sigma}_t\mathbf{z}_t=O_p(p)$. 
We then get
\begin{equation}\label{trace_difference_R}
\tr\Big((\widetilde{\mathbf{R}}_t-\check{{\mathbf{R}}}_t)(\widetilde{\mathbf{R}}_t-\check{{\mathbf{R}}}_t)^T\Big)=o_p(p).
\end{equation}
In addition, by $\mathbf{P}_t\geq \delta \mathbf{I}$, and $\check{\mathbf{P}}_t\geq \delta \mathbf{I}$, we have
$$
\aligned
\tr\Big((\widetilde{\mathbf{R}}_t-\check{{\mathbf{R}}}_t)(\widetilde{\mathbf{R}}_t-\check{{\mathbf{R}}}_t)^T/p\Big)\leq&2\tr\Big(\widetilde{\mathbf{R}}_t\widetilde{\mathbf{R}}_t^T/p+\check{{\mathbf{R}}}_t\check{{\mathbf{R}}}_t^T/p\Big)\leq \frac{4}{\delta}\mathbf{R}_{t}^T\mathbf{R}_t/p.
\endaligned
$$
By the independence between $(\mathbf{z}_t)$ and $(\boldsymbol{\Sigma}_t)$ and Assumption \ref{asump1}(iii), we have $E(\mathbf{R}_{t}^T\mathbf{R}_t/p)=\tr(\overline{\boldsymbol{\Sigma}})/p=O(1)$. 
By \eqref{trace_difference_R} and the dominated convergence theorem, 
$$ E\Bigg(\tr\Big((\widetilde{\mathbf{R}}_t-\check{{\mathbf{R}}}_t)(\widetilde{\mathbf{R}}_t-\check{{\mathbf{R}}}_t)^T\Big)/p\Bigg)=o(1).
$$
It follows that
$$
\frac{1}{pn} E\Bigg(\tr\Big((\widetilde{\mathbf{R}}-\check{{\mathbf{R}}})(\widetilde{\mathbf{R}}-\check{{\mathbf{R}}})^T\Big)\Bigg)= E\Bigg(\tr\Big((\widetilde{\mathbf{R}}_t-\check{{\mathbf{R}}}_t)(\widetilde{\mathbf{R}}_t-\check{{\mathbf{R}}}_t)^T\Big)/p\Bigg)=o(1).
$$
By Markov's inequality,
we get
\begin{equation}\label{part2_check}
\frac{1}{pn}\tr\Big((\widetilde{\mathbf{R}}-\check{{\mathbf{R}}})(\widetilde{\mathbf{R}}-\check{{\mathbf{R}}})^T\Big)=o_p(1).
\end{equation}
By \eqref{L_check_diff}, \eqref{trS_ntilde} and \eqref{part2_check}, we have
\begin{equation}\label{L_replace1}
L(F^{\widetilde{\mathbf{S}}_n},F^{\check{{\mathbf{S}}}_n})=o_p(1).
\end{equation}

\underline{\emph{Step Two}}: We denote by $\mathcal{F}_t$ the $\sigma$-algebra generated by $\{\mathbf{z}_s,\infty<s\leq t\}$. For a $p\times p$ matrix ${\mathcal{O}}_{t}$ to be determined, which satisfies that
\begin{equation}\label{def_o}
{\mathcal{O}}_{t} \text{ is }\mathcal{F}_{t-1} \text{-measurable},\quad\text{and } \;{\mathcal{O}}_{t}{\mathcal{O}}_{t}^T={\mathcal{O}}_{t}^T{\mathcal{O}}_{t}=\mathbf{I},
\end{equation}
we perform orthogonal transformation on $\mathbf{z}_t$ and get $\boldsymbol{\zeta}_t={\mathcal{O}}_{t}\mathbf{z}_t$.
We then define 
\begin{equation}\label{R0*_def}
\boldsymbol{\mathcal{R}}^0_{t}=\overline{\boldsymbol{\Sigma}}^{1/2}\boldsymbol{\zeta}_t,\quad \boldsymbol{\mathcal{R}}^0=(\boldsymbol{\mathcal{R}}^0_{1}, ..., \boldsymbol{\mathcal{R}}^0_{n}), \quad\text{and } \;\boldsymbol{\mathcal{S}}_n^{0}=\frac{1}{n}\sum_{t=1}^n\boldsymbol{\mathcal{R}}^0_{t}(\boldsymbol{\mathcal{R}}^0_{t})^T.
\end{equation}
By \eqref{def_o} and the assumption that $\mathbf{z}_{t}\underset{\text{i.i.d.}}\sim  {N}(0, \mathbf{I})$, we have 
\begin{equation}\label{invariant_normal}
\boldsymbol{\zeta}_t\underset{\text{i.i.d.}}\sim  {N}(0, \mathbf{I}).
\end{equation}
By Theorem~1 of \cite{MP67}, $F^{\mathbf{S}^0_n}\overset{\text{P}}\to F$, and 
$F^{\boldsymbol{\mathcal{S}}_n^{0}}\overset{\text{P}}\to F$. Hence,
\begin{equation}\label{S_n^*0vsS_n^0}
L(F^{{\mathbf{S}}^0_n}, F^{\boldsymbol{\mathcal{S}}_n^{0}})=o_p(1).
\end{equation}
By \eqref{L_replace1}, \eqref{S_n^*0vsS_n^0} and the triangle inequality, to show Theorem~\ref{recovered}, it suffices to show that
\begin{equation}\label{L_S*0_check}
L(F^{\boldsymbol{\mathcal{S}}_n^{0}},F^{\check{{\mathbf{S}}}_n})=o_p(1).
\end{equation}
By Corollary A.42 of \cite{BS10} again, 
\begin{equation}\label{L4_recovered}
L^4(F^{\check{{\mathbf{S}}}_n},  F^{\boldsymbol{\mathcal{S}}^{0}_n})\leq \frac{2}{p}\tr (\check{{\mathbf{S}}}_n+\boldsymbol{\mathcal{S}}^{0}_n)\cdot \frac{1}{pn}\tr\Big((\check{{\mathbf{R}}}-\boldsymbol{\mathcal{R}}^0)(\check{{\mathbf{R}}}-\boldsymbol{\mathcal{R}}^0)^T\Big). 
\end{equation}
We have $E\big(\tr(\boldsymbol{\mathcal{S}}_n^{0}/p)\big)=\tr(\overline{\boldsymbol{\Sigma}})/p=O(1)$, hence 
\begin{equation}\label{tr_s0}
\tr(\boldsymbol{\mathcal{S}}_n^{0}/p)=O_p(1).
\end{equation}
Combining \eqref{trS_ntilde} and \eqref{tr_s0} yields
\begin{equation}\label{tr_tildS}
\tr(\check{{\mathbf{S}}}_n/p+\boldsymbol{\mathcal{S}}_n^{0}/p)=O_{p}(1).
\end{equation}
Define
\begin{equation}\label{tilde_SigmaL}
\mathbf{Q}_{t}=\check{\mathbf{P}}_{t}^{-1/2}\boldsymbol{\Sigma}_t^{1/2}.
\end{equation} We have $\check{{\mathbf{R}}}_t=\mathbf{Q}_{t}\mathbf{z}_t$, and $\check{{\boldsymbol{\Sigma}}}_t=\mathbf{Q}_{t}(\mathbf{Q}_{t})^T$. 
We will show that for some ${\mathcal{O}}_{t}$ satisfying
\eqref{def_o},
\begin{equation}\label{E_tildR}
\frac{1}{p}E\Bigg(\tr\Big((\mathbf{Q}_{t}{\mathcal{O}}_{t}^T-\overline{\boldsymbol{\Sigma}}^{1/2})(\mathbf{Q}_{t}{\mathcal{O}}_{t}^T-\overline{\boldsymbol{\Sigma}}^{1/2})^T\Big)\Bigg)=o(1).
\end{equation}
Then by the facts that $\check{{\boldsymbol{\Sigma}}}_{t}$ and ${\mathcal{O}}_{t}$ are $\mathcal{F}_{t-1}$-measurable, we have
\begin{equation}\label{Jesen_eq_tilde}
\aligned
&\frac{1}{np}E\Bigg(\tr\Big((\check{{\mathbf{R}}}-\boldsymbol{\mathcal{R}}^0)(\check{{\mathbf{R}}}-\boldsymbol{\mathcal{R}}^0)^T\Big)\Bigg)\\
=&\frac{1}{p}E\Bigg(\tr\Big((\check{{\mathbf{R}}}_t-\boldsymbol{\mathcal{R}}^0_{t})(\check{{\mathbf{R}}}_t-\boldsymbol{\mathcal{R}}^0_{t})^T\Big)\Bigg)\\
=&\frac{1}{p}E\Bigg(\tr\Big((\mathbf{Q}_{t}-\overline{\boldsymbol{\Sigma}}^{1/2}{\mathcal{O}}_{t})\mathbf{z}_t\mathbf{z}_t^T(\mathbf{Q}_{t}-\overline{\boldsymbol{\Sigma}}^{1/2}{\mathcal{O}}_{t})^T\Big)\Bigg)\\
=&\frac{1}{p}E\Bigg(\tr\Big((\mathbf{Q}_{t}{\mathcal{O}}_{t}^T-\overline{\boldsymbol{\Sigma}}^{1/2})(\mathbf{Q}_{t}{\mathcal{O}}_{t}^T-\overline{\boldsymbol{\Sigma}}^{1/2})^T\Big)\Bigg)=o(1),
\endaligned
\end{equation}
which implies that
\begin{equation}\label{trace_diff}
\frac{1}{pn}\tr\Big((\check{{\mathbf{R}}}-\boldsymbol{\mathcal{R}}^0)(\check{{\mathbf{R}}}-\boldsymbol{\mathcal{R}}^0)^T\Big)=o_p(1).
\end{equation}
The desired bound \eqref{L_S*0_check} then follows from \eqref{L4_recovered}, \eqref{tr_tildS} and \eqref{trace_diff}. 

\underline{\emph{Step Three}}: It remains to show that there exists $\mathcal{O}_{t}$ satisfying \eqref{def_o} and \eqref{E_tildR}.
Because $M_p=o(\sqrt{p})\ll p$, with probability one,  for all $p$ large enough,\\
$\text{rank}(\sum_{j=1}^{M_p} {b}^{j-1} \mathbf{R}_{t-j}\mathbf{R}_{t-j}^T)= M_p$. Write
$$
\aligned
&\sum_{j=1}^{M_p} {b}^{j-1} \mathbf{R}_{t-j}\mathbf{R}_{t-j}^T=:\mathbf{U}\mathbf{\Lambda}\mathbf{U}^T,
\endaligned
$$
where $\mathbf{\Lambda}=\diag(\lambda_1, ..., \lambda_{M_p})$ and $\mathbf{U}=(\mathbf{u}_1, ... , \mathbf{u}_{M_p})$  are the  nonzero eigenvalues and the corresponding eigenvectors of $\sum_{j=1}^{M_p} {b}^{j-1} \mathbf{R}_{t-j}\mathbf{R}_{t-j}^T$, respectively. 
Recall that 
$$\aligned
\check{\mathbf{P}}_{t}=&\frac{1-{a}-{b}+{a}{b}^{M_p}}{1-{b}}\mathbf{I}+\sum_{j=1}^{M_p} {a} {b}^{j-1} \mathbf{R}_{t-j}\mathbf{R}_{t-j}^T=\frac{1-{a}-{b}+{a}{b}^{M_p}}{1-{b}}\mathbf{I}+{a}\mathbf{U}\boldsymbol{\Lambda}\mathbf{U}^T.
\endaligned
$$
We have
\begin{equation}\label{Woodbury}
\aligned
\check{\mathbf{P}}_t^{-1/2}=&\mathbf{U}\Big(a\boldsymbol{\Lambda}+\frac{1-{a}-{b}+{a}{b}^{M_p}}{1-{b}}\mathbf{I}\Big)^{-1/2}\mathbf{U}^T\\
&+\sqrt{\frac{1-{b}}{1-{a}-{b}+{a}{b}^{M_p}}}(\mathbf{I}-\mathbf{U}\mathbf{U}^T).\endaligned
\end{equation}
By \eqref{BEKK}, we have
\begin{equation}\label{sigma_inf_sum_presentation}
\boldsymbol{\Sigma}_{t}=\frac{1-a-b}{1-b}\overline{\boldsymbol{\Sigma}}+\sum_{s=1}^\infty a b^{s-1}\mathbf{R}_{t-s}\mathbf{R}_{t-s}^T.
\end{equation}
Hence
\begin{equation}\label{Sig_decomp}
\aligned
&\boldsymbol{\Sigma}_{t}\\
=&\Bigg(\frac{1-a-b+ab^{M_p}}{1-{b}}\overline{\boldsymbol{\Sigma}}+a \mathbf{U}\boldsymbol{\Lambda}\mathbf{U}^T\Bigg)+\Bigg(\sum_{j=M_p+1}^\infty a b^{j-1} \mathbf{R}_{t-j}\mathbf{R}_{t-j}^T-\frac{ab^{M_p}}{1-{b}}\overline{\boldsymbol{\Sigma}}\Bigg)
\\
=&:I_{t}+II_t.
\endaligned
\end{equation}
Recall that $\check{{\boldsymbol{\Sigma}}}_t=\check{\mathbf{P}}_t^{-1/2}\boldsymbol{\Sigma}_t\check{\mathbf{P}}_t^{-1/2}$. We have
\begin{equation}\label{tildeSig_decomp}
\aligned
\check{{\boldsymbol{\Sigma}}}_t
=&\check{\mathbf{P}}_t^{-1/2}I_t\check{\mathbf{P}}_t^{-1/2}+\check{\mathbf{P}}_t^{-1/2}II_t\check{\mathbf{P}}_t^{-1/2}
.\endaligned
\end{equation}
By \eqref{Woodbury}, 
\begin{equation}\label{I_t_decom}
\aligned
&\check{\mathbf{P}}_t^{-1/2}I_t\check{\mathbf{P}}_t^{-1/2}=\overline{\boldsymbol{\Sigma}}-\overline{\boldsymbol{\Sigma}}\mathbf{U}\mathbf{U}^T-\mathbf{U}\mathbf{U}^T\overline{\boldsymbol{\Sigma}}+\mathbf{U}\mathbf{U}^T\overline{\boldsymbol{\Sigma}}\mathbf{U}\mathbf{U}^T+\mathcal{E}_t
,\endaligned
\end{equation}
where
$$
\aligned&\mathcal{E}_t\\
=&\frac{1-{a}-{b}+{a}{b}^{M_p}}{1-{b}}\mathbf{U}\Big(a\boldsymbol{\Lambda}+\frac{1-{a}-{b}+{a}{b}^{M_p}}{1-b}\mathbf{I}\Big)^{-1/2}\mathbf{U}^T\overline{\boldsymbol{\Sigma}}\mathbf{U}\Big(a\boldsymbol{\Lambda}+\frac{1-{a}-{b}+{a}{b}^{M_p}}{1-b}\mathbf{I}\Big)^{-1/2}\mathbf{U}^T\\
&+ \sqrt{\frac{1-{a}-{b}+{a}{b}^{M_p}}{1-b}}
\Big(\mathbf{I}-\mathbf{U}\mathbf{U}^T\Big)\overline{\boldsymbol{\Sigma}}\mathbf{U}\Big(a\boldsymbol{\Lambda}+\frac{1-{a}-{b}+{a}{b}^{M_p}}{1-b}\mathbf{I}\Big)^{-1/2}\mathbf{U}^T\\
&+\sqrt{\frac{1-{a}-{b}+{a}{b}^{M_p}}{1-b}}
\mathbf{U}\Big(a\boldsymbol{\Lambda}+\frac{1-{a}-{b}+{a}{b}^{M_p}}{1-b}\mathbf{I}\Big)^{-1/2}\mathbf{U}^T\overline{\boldsymbol{\Sigma}}\Big(\mathbf{I}-\mathbf{U}\mathbf{U}^T\Big)\\
&+
a\mathbf{U}\Big(a\boldsymbol{\Lambda}+\frac{1-{a}-{b}+{a}{b}^{M_p}}{1-b}\mathbf{I}\Big)^{-1/2}\boldsymbol{\Lambda}\Big(a\boldsymbol{\Lambda}+\frac{1-{a}-{b}+{a}{b}^{M_p}}{1-b}\mathbf{I}\Big)^{-1/2}\mathbf{U}^T\\
=:&\mathcal{E}_{1t}+\mathcal{E}_{2t}+\mathcal{E}_{3t}+\mathcal{E}_{4t}.
\endaligned
$$
Because $\mathbf{U}^T\mathbf{U}=\mathbf{I}$, we have $$
\tr(\mathcal{E}_{2t})=\tr(\mathcal{E}_{3t})=0.$$ Moreover, because $M_p\ll p$ and $\|\overline{\boldsymbol{\Sigma}}\|=O(1)$, 
\begin{equation}\label{tr_usigmau}
0\leq\tr(\mathbf{U}^T\overline{\boldsymbol{\Sigma}}\mathbf{U})=\sum_{i=1}^{M_p}\mathbf{u}_i^T\overline{\boldsymbol{\Sigma}}\mathbf{u}_i\leq M_p\|\overline{\boldsymbol{\Sigma}}\|=o(p).
\end{equation}
Furthermore, by \eqref{QMLE_bound}, $M_p\ll p$ and the fact that if $\mathbf{A}\geq 0,\mathbf{B}\geq 0$, then $\tr(\mathbf{A}\mathbf{B})\geq 0$, we have 
$$
\aligned
0\leq&\tr(\mathcal{E}_{1t})\leq 
\tr(\mathbf{U}^T\overline{\boldsymbol{\Sigma}}\mathbf{U})=o(p),\\
0\leq&\tr(\mathcal{E}_{4t})=\tr\Bigg(a\boldsymbol{\Lambda}\Big(a\boldsymbol{\Lambda}+\frac{1-{a}-{b}+{a}{b}^{M_p}}{1-b}\mathbf{I}\Big)^{-1}\Bigg)\leq M_p=o(p).
\endaligned
$$
Combining the results above yields
$$\tr(\mathcal{E}_t/p)=o(1).$$ 
Moreover, by \eqref{tr_usigmau} and that $\mathbf{U}^T\mathbf{U}=\mathbf{I}$, we have
$$
\aligned
0\leq &\tr(\mathbf{U}\mathbf{U}^T\overline{\boldsymbol{\Sigma}}\mathbf{U}\mathbf{U}^T)=\tr(\mathbf{U}\mathbf{U}^T\overline{\boldsymbol{\Sigma}})=\tr(\overline{\boldsymbol{\Sigma}}\mathbf{U}\mathbf{U}^T)=\tr(\mathbf{U}^T\overline{\boldsymbol{\Sigma}}\mathbf{U})=o(p).
\endaligned$$
Plugging the estimates above  into \eqref{I_t_decom} yields
\begin{equation}\label{tildeSig_decomp_1}
\tr(\check{\mathbf{P}}_t^{-1/2}I_t\check{\mathbf{P}}_t^{-1/2})-\tr(\overline{\boldsymbol{\Sigma}})=o(p).
\end{equation}

About term $\check{\mathbf{P}}_t^{-1/2}II_t\check{\mathbf{P}}_t^{-1/2}$, because $\check{\mathbf{P}}_t\geq (1-a-b)\mathbf{I}$, we have
$$\aligned
&\Big|\tr\Big(\check{\mathbf{P}}_t^{-1/2}II_t\check{\mathbf{P}}_t^{-1/2}\Big)\Big|\\
\leq &\frac{ab^{M_p}}{(1-a-b)(1-b)}\tr(\overline{\boldsymbol{\Sigma}})+\frac{b^{M_p}}{1-a-b} \tr\Big(\sum_{s=1}^{\infty}ab^{s-1} \mathbf{R}_{t-M_p-s}\mathbf{R}_{t-M_p-s}^T\Big)\\
\leq& \frac{ab^{M_p}}{(1-a-b)(1-{b})}\tr(\overline{\boldsymbol{\Sigma}})+\frac{b^{M_p}}{1-a-b}\tr(\boldsymbol{\Sigma}_{t-M_p}).
\endaligned
$$
We have $E\big(\tr(\boldsymbol{\Sigma}_{t-M_p})\big)=\tr\big(\overline{\boldsymbol{\Sigma}}\big)=O(p)$, hence $\tr(\boldsymbol{\Sigma}_{t-M_p})=O_p(p)$.
By the assumptions that $M_p\to\infty$ and $a+b<1-\delta$, we get
\begin{equation}\label{tildeSig_decomp_4}
\aligned
\Big|\tr\Big(\check{\mathbf{P}}_t^{-1/2}II_t\check{\mathbf{P}}_t^{-1/2}\Big)\Big|=&O(b^{M_p})\Big(\tr(\overline{\boldsymbol{\Sigma}})+\tr(\boldsymbol{\Sigma}_{t-M_p})\Big)\\
=&O(b^{M_p})O_p(p)=o_p(p).
\endaligned
\end{equation}
Combining \eqref{tildeSig_decomp_1} and \eqref{tildeSig_decomp_4} yields
\begin{equation}\label{dev_large2}
\Big|\tr(\check{{\boldsymbol{\Sigma}}}_t)/p-\tr(\overline{\boldsymbol{\Sigma}})/p\Big|=o_p(1). 
\end{equation}

We now define $\mathcal{O}_t$ that satisfies \eqref{def_o} and \eqref{E_tildR}. Let
$$\mathbf{G}_t=\sqrt{(1-a-b)/(1-b)}\check{\mathbf{P}}_t^{-1/2}\overline{\boldsymbol{\Sigma}}^{1/2}.$$ We have $\mathbf{G}_t\mathbf{G}_t^T=\big((1-a-b)/(1-b)\big)\check{\mathbf{P}}_t^{-1/2}\overline{\boldsymbol{\Sigma}}\check{\mathbf{P}}_t^{-1/2}.$  By \eqref{sigma_inf_sum_presentation}, $\boldsymbol{\Sigma}_{t}\geq \big((1-a-b)/(1-b)\big)\overline{\boldsymbol{\Sigma}}$.
Hence
\begin{equation}\label{Sigma_lowerbound}
\check{{\boldsymbol{\Sigma}}}_t\geq \mathbf{G}_t\mathbf{G}_t^T.
\end{equation}
Define $$ \mathcal{Q}_t=\mathbf{G}_t \Big(\mathbf{I}+\mathbf{G}_t^{-1}(\check{{\boldsymbol{\Sigma}}}_t-\mathbf{G}_t\mathbf{G}_t^T)(\mathbf{G}_t^T)^{-1}\Big)^{1/2},
$$
and 
$${\mathcal{O}}_{t}=\mathcal{Q}_t^T\big(\mathbf{Q}_t^T\big)^{-1},
$$
where, recall that, $\mathbf{Q}_t$ is defined in \eqref{tilde_SigmaL}. By definition, ${\mathcal{O}}_{t}$ is $\mathcal{F}_{t-1}$-measurable. Moreover, it is straightforward to verify that  $
\mathcal{Q}_t\mathcal{Q}_t^T=
\mathbf{Q}_t\mathbf{Q}_t^T=\check{{\boldsymbol{\Sigma}}}_t$, from which we get that ${\mathcal{O}}_{t}{\mathcal{O}}^T_{t}={\mathcal{O}}^T_{t}{\mathcal{O}}_{t}=\mathbf{I}$. Therefore, ${\mathcal{O}}_{t}$ satisfies \eqref{def_o}. 

It remains to show \eqref{E_tildR}.
By \eqref{Woodbury},  
\begin{equation}\label{P_lowerbound}
\check{\mathbf{P}}_{t}^{-1/2}\geq \sqrt{\frac{1-b}{1-{a}-{b}+{a}{b}^{M_p}}} \Big(\mathbf{I}-\mathbf{U}\mathbf{U}^T\Big).
\end{equation}
By \eqref{Sigma_lowerbound} and \eqref{P_lowerbound}, 
\begin{equation}\label{tr_1/2}
\aligned
&\tr\Big(\mathcal{Q}_t^T\overline{\boldsymbol{\Sigma}}^{1/2}\Big)=\tr\Big(\mathcal{Q}_t\overline{\boldsymbol{\Sigma}}^{1/2}\Big)\\
=&\sqrt{\frac{1-a-b}{1-b}}\tr\Big(\big(\mathbf{I}+\mathbf{G}_t^{-1}(\check{{\boldsymbol{\Sigma}}}_t-\mathbf{G}_t\mathbf{G}_t^T)(\mathbf{G}_t^T)^{-1}\big)^{1/2}\overline{\boldsymbol{\Sigma}}^{1/2}\check{\mathbf{P}}_t^{-1/2}\overline{\boldsymbol{\Sigma}}^{1/2}\Big)\\
\geq& \sqrt{\frac{1-a-b}{1-b}}\tr \big(\check{\mathbf{P}}_t^{-1/2}\overline{\boldsymbol{\Sigma}}\big)\\
\geq& \sqrt{\frac{1-{a}-{b}}{1-{a}-{b}+{a}{b}^{M_p}}} \tr \Big(\big(\mathbf{I}-\mathbf{U}\mathbf{U}^T\big)\overline{\boldsymbol{\Sigma}}\Big)\\
=&\tr(\overline{\boldsymbol{\Sigma}})+o(p),
\endaligned
\end{equation}
where the last equation holds by \eqref{tr_usigmau} and the assumptions that $M_p\to\infty$ and $a+b<1-\delta$.
By the definition of $\mathcal{O}_t$, \eqref{dev_large2} and \eqref{tr_1/2},
we get that
$$
\aligned
0\leq &\frac{1}{p}\tr\Big((\mathbf{Q}_{t}{\mathcal{O}}_{t}^T-\overline{\boldsymbol{\Sigma}}^{1/2})(\mathbf{Q}_{t}{\mathcal{O}}_{t}^T-\overline{\boldsymbol{\Sigma}}^{1/2})^T\Big)
=\frac{1}{p}\Bigg(\tr\left(\Big(\mathbf{Q}_t-\overline{\boldsymbol{\Sigma}}^{1/2}{\mathcal{O}}_{t}\Big)\Big(\mathbf{Q}_t-\overline{\boldsymbol{\Sigma}}^{1/2}{\mathcal{O}}_{t}\Big)^T\right)\Bigg)\\
=&\frac{1}{p}\Bigg(\tr\left(\Big(\mathcal{Q}_t-\overline{\boldsymbol{\Sigma}}^{1/2}\Big)\Big(\mathcal{Q}_t-\overline{\boldsymbol{\Sigma}}^{1/2}\Big)^T\right)\Bigg)\\
=&\Big(\tr(\check{{\boldsymbol{\Sigma}}}_{t})/p-\tr(\overline{\boldsymbol{\Sigma}})/p\Big)+\Big(\tr(\overline{\boldsymbol{\Sigma}})/p-\tr(\mathcal{Q}_{t}\overline{\boldsymbol{\Sigma}}^{1/2})/p\Big)+\Big(\tr(\overline{\boldsymbol{\Sigma}})/p-\tr\big(\mathcal{Q}_t^T\overline{\boldsymbol{\Sigma}}^{1/2}\big)/p\Big)\\
=&o_p(1).
\endaligned
$$
In addition, by \eqref{bound_sigma_ii}, we have
$$
\aligned
&\tr\Big(\big(\mathcal{Q}_t-\overline{\boldsymbol{\Sigma}}^{1/2}\big)\big(\mathcal{Q}_t-\overline{\boldsymbol{\Sigma}}^{1/2}\big)^T\Big)/p\\
=&\sum_{1\leq i,j\leq p}\Big((\mathcal{Q}_t)_{ij}-(\overline{\boldsymbol{\Sigma}}^{1/2})_{ij}\Big)^2/p\\
\leq& 2\sum_{1\leq i,j\leq p}\Big((\mathcal{Q}_t)_{ij}\Big)^2/p+2\sum_{1\leq i,j\leq p}\Big((\overline{\boldsymbol{\Sigma}}^{1/2})_{ij}\Big)^2/p\\
= &2 \tr(\check{{\boldsymbol{\Sigma}}}_t/p)+2\tr(\overline{\boldsymbol{\Sigma}}/p)\\
\leq& \frac{2}{\delta} \Big(\tr(\boldsymbol{\Sigma}_{t})/p+\tr(\overline{\boldsymbol{\Sigma}}/p)\Big).
\endaligned
$$ 
By Assumption \ref{asump1}(iii), we have
$E\big(\tr(\boldsymbol{\Sigma}_{t})/p\big)+\tr(\overline{\boldsymbol{\Sigma}}/p)=2\tr(\overline{\boldsymbol{\Sigma}})/p\leq 2\|\overline{\boldsymbol{\Sigma}}\|<2C.$ By the dominated convergence theorem again,
the bound \eqref{E_tildR} follows.
$\qed$


%
%

\section*{Funding}
Research is supported in part by RGC grants GRF~16304019, GRF~15302321 and
GRF~16304521  of the HKSAR.

\section*{Supplement Materials}
{\bf Supplement to ``High-dimensional covariance matrices under dynamic volatility models: asymptotics and shrinkage estimation''}.

{This supplement contains the proofs of Theorems \ref{reducable}, \ref{nonreducable}, \ref{thetag} and \ref{NLS_ADJ_cov} and  Corollary \ref{NLS_ADJ}.
}



}}}

\bibliographystyle{apalike}
\bibliography{DCC_RMT_aos-template}       



\newpage
\setcounter{page}{1}
\clearpage

\begin{center}
    {\Large \bf Supplement to ``High-dimensional covariance matrices under dynamic volatility models: asymptotics and shrinkage estimation''}
    
    \bigskip
  
  {\large Yi Ding and Xinghua Zheng }

\end{center}

\pagenumbering{arabic}
\numberwithin{equation}{section}

\setcounter{equation}{0}
\renewcommand{\theequation}{A.\arabic{equation}}
For ease of notation, we drop the subscript~$p$ in~$a_p$ and~$b_p$ in the following proofs. 

\section*{Proof of~Theorem~\ref{reducable}}

Recall that $\eta(a,b,p)=\big(a/(1-a-b)\big)\min(\sqrt{p(1-a-b)},1)$. 
Under the assumption $\eta(a,b,p)\to 0$, we separately consider the following two cases:
$$
\begin{cases}
\text{Case I: when }& a/(1-a-b)\to 0;\\
\text{Case II: when }&a/(1-a-b)>c\text{ and }a^2p/(1-a-b)\to 0.
\end{cases}
$$

We first consider Case I.
By \eqref{sigma_inf_sum_presentation} and the L$\ddot{\text{o}}$wner-Heinz inequality, \begin{equation}\label{ineq_sigmat_half}
\boldsymbol{\Sigma}^{1/2}_{t}\geq \sqrt{1-a/(1-b)}\overline{\boldsymbol{\Sigma}}^{1/2}.
\end{equation}
By Corollary A.42 of \cite{BS10}, 
\begin{equation}\label{L4}
L^4(F^{\mathbf{S}_n},  F^{\mathbf{S}^{0}_n})\leq \frac{2}{p^2n}\tr (\mathbf{S}_n+\mathbf{S}^{0}_n)\cdot \tr\Big((\mathbf{R}-\mathbf{R}^0)(\mathbf{R}-\mathbf{R}^0)^T\Big), 
\end{equation}
where~$\mathbf{R}=(\mathbf{R}_1, ..., \mathbf{R}_n)$ and~$\mathbf{R}^0=(\mathbf{R}^0_{1}, ..., \mathbf{R}^0_{ n})$. Under model \eqref{BEKK}, by the independence between $(\boldsymbol{\Sigma}_t)$ and $(\mathbf{z}_t)$, we have $$E\Big(\tr (\mathbf{S}_n)/p\Big)=E\Big(\tr(\mathbf{S}^0_n)/p\Big)=\tr(\overline{\boldsymbol{\Sigma}})/p.$$ Therefore, by Assumption \ref{asump1}(iii), we have
\begin{equation}\label{part1}
\frac{1}{p}\tr (\mathbf{S}_n+\mathbf{S}^{0}_n)=O_{p}(1). 
\end{equation}
By \eqref{ineq_sigmat_half} and the fact that if $\mathbf{A}\geq 0$ and $\mathbf{B}\geq 0$, then $\tr(\mathbf{A}\mathbf{B})\geq 0$, we get
$$
\aligned
&\tr\Big(\boldsymbol{\Sigma}_{t}^{1/2}\overline{\boldsymbol{\Sigma}}^{1/2}\Big)\geq &\sqrt{1-a/(1-b)}\tr(\overline{\boldsymbol{\Sigma}}).
\endaligned
$$
Under model \eqref{BEKK}, we have $E\big(\tr(\boldsymbol{\Sigma}_{t})\big)=\tr(\overline{\boldsymbol{\Sigma}})$. Combining these results and using the definitions of $\mathbf{R}_t$ and $\mathbf{R}^0_{t}$, we get
$$
\aligned
0\leq &E\Bigg(\frac{1}{np}\tr\Big((\mathbf{R}-\mathbf{R}^0)(\mathbf{R}-\mathbf{R}^0)^T\Big)\Bigg)\\
=&\frac{1}{p}E\Bigg(\tr\Big((\boldsymbol{\Sigma}_t^{1/2}-\overline{\boldsymbol{\Sigma}}^{1/2})(\boldsymbol{\Sigma}_t^{1/2}-\overline{\boldsymbol{\Sigma}}^{1/2})^T\Big)\Bigg)\\
=&\frac{2}{p}\Bigg(\tr(\overline{\boldsymbol{\Sigma}})-E\tr\Big(\boldsymbol{\Sigma}_{t}^{1/2}\overline{\boldsymbol{\Sigma}}^{1/2}\Big)\Bigg)\\
\leq &2\Big(1-\sqrt{1-a/(1-b)}\Big)\tr(\overline{\boldsymbol{\Sigma}}/p)=o(1),
\endaligned
$$
where the last equation holds by the assumption that $a/(1-a-b)\to 0$, which implies that $a/(1-b)\to 0$ as $p\to \infty$. 
By Markov's inequality,
\begin{equation}\label{part2_1}
\frac{1}{np}\tr\Big((\mathbf{R}-\mathbf{R}^0)(\mathbf{R}-\mathbf{R}^0)^T\Big)=o_{p}(1).
\end{equation}
The desired bound \eqref{thm1_1}  follows from \eqref{L4}, \eqref{part1} and \eqref{part2_1}. 

Next, we consider Case II. 
Define $$\mathbf{W}_t=\overline{\boldsymbol{\Sigma}}^{1/2}\Bigg((1-a-b)/(1-b)\mathbf{I}+\sum_{s=1}^\infty a b^{s-1}\overline{\boldsymbol{\Sigma}}^{-1/2}\mathbf{R}_{t-s}\mathbf{R}_{t-s}^T\overline{\boldsymbol{\Sigma}}^{-1/2}\Bigg)^{1/2},$$
and $${\mathfrak{O}}_{t}=\mathbf{W}_t^{-1}\boldsymbol{\Sigma}^{1/2}_t.$$
It is straightforward to verify that $\mathbf{W}_t\mathbf{W}_t^T=\boldsymbol{\Sigma}_t$, from which we get that
 \begin{equation}\label{def_ow}
 {\mathfrak{O}}_{t}({\mathfrak{O}}_{t})^T=({\mathfrak{O}}_{t})^T{\mathfrak{O}}_{t}=\mathbf{I}.
 \end{equation}
Define 
\begin{equation}\label{R0w_def}
\boldsymbol{\varepsilon}_t={\mathfrak{O}}_{t}\mathbf{z}_t, \quad \mathfrak{R}_{t}^{0}=\overline{\boldsymbol{\Sigma}}^{1/2}\boldsymbol{\varepsilon}_t,\quad \mathfrak{R}^{0}=(\mathfrak{R}_1^{0}, ..., \mathfrak{R}_{n}^{0}), \quad\text{and }\, \mathfrak{S}_n^{0}=\frac{1}{n}\sum_{t=1}^n\mathfrak{R}^{0}_{t}(\mathfrak{R}^{0}_{t})^T.
\end{equation}
Note that $\mathfrak{O}_t$ is $\mathcal{F}_{t-1}$-measurable. By \eqref{def_ow} and the assumption that $\mathbf{z}_{t}\underset{\text{i.i.d.}}\sim  {N}(0, \mathbf{I})$, we have
$$
\boldsymbol{\varepsilon}_t\underset{\text{i.i.d.}}\sim  {N}(0, \mathbf{I}).
$$
It follows from Theorem~1 of \cite{MP67} that $F^{\mathbf{S}^0_n}\overset{\text{P}}\to F$, and 
$F^{\mathfrak{S}_n^{0}}\overset{\text{P}}\to F$. Hence,
\begin{equation}\label{S_n^w0vsS_n^0}
L(F^{{\mathbf{S}}^0_n}, F^{\mathfrak{S}_n^{0}})=o_p(1).
\end{equation}
By \eqref{S_n^w0vsS_n^0} and the triangle inequality, to show \eqref{thm1_1}, it suffices to show that
\begin{equation}\label{L_Sw0_check}
L(F^{\mathfrak{S}_n^{0}},F^{{{\mathbf{S}}}_n})=o_p(1).
\end{equation}
By Corollary A.42 of \cite{BS10}, 
\begin{equation}\label{L4_recovered_w}
L^4(F^{{{\mathbf{S}}}_n},  F^{\mathfrak{S}^{0}_n})\leq \frac{2}{p}\tr ({{\mathbf{S}}}_n+\mathfrak{S}^{0}_n)\cdot \frac{1}{pn}\tr\Big(({{\mathbf{R}}}-\mathfrak{R}^{0})({{\mathbf{R}}}-\mathfrak{R}^{0})^T\Big). 
\end{equation}
We have $E\big(\tr(\mathfrak{S}_n^{0}/p)\big)=\tr(\overline{\boldsymbol{\Sigma}})/p=O(1)$. By Markov's inequality, 
\begin{equation}\label{tr_s0w}
\tr(\mathfrak{S}_n^{0}/p)=O_p(1).
\end{equation}
By the independence between $(\boldsymbol{\Sigma}_t)$ and $(\mathbf{z}_t)$, $E\big(\tr(\mathbf{S}_{n})\big)=E\big(\tr(\boldsymbol{\Sigma}_{t})\big)=\tr(\overline{\boldsymbol{\Sigma}})$. Therefore, \eqref{part1} holds. 
Combining \eqref{part1} and \eqref{tr_s0w} yields
\begin{equation}\label{tr_tildS_w}
\tr({{\mathbf{S}}}_n/p+\mathfrak{S}_n^{0}/p)=O_{p}(1).
\end{equation}
We will show that
\begin{equation}\label{E_tildR_w}
\frac{1}{p}E\Bigg(\tr\Big(({\mathbf{W}}_t-\overline{\boldsymbol{\Sigma}}^{1/2})({\mathbf{W}}_t-\overline{\boldsymbol{\Sigma}}^{1/2})^T\Big)\Bigg)=o(1).
\end{equation}
Then by the facts that ${{\boldsymbol{\Sigma}}}_{t}$ and ${\mathfrak{O}}_{t}$ are $\mathcal{F}_{t-1}$-measurable, we have
\begin{equation}\label{Jesen_eq_tilde_W}
\aligned
&\frac{1}{np}E\Bigg(\tr\Big(({{\mathbf{R}}}-\mathfrak{R}^{0})({{\mathbf{R}}}-\mathfrak{R}^{0})^T\Big)\Bigg)\\
=&\frac{1}{p}E\Bigg(\tr\Big(({{\mathbf{R}}}_t-\mathfrak{R}^{0}_{t})({{\mathbf{R}}}_t-\mathfrak{R}^{0}_{t})^T\Big)\Bigg)\\
=&\frac{1}{p}E\Bigg(\tr\Big((\boldsymbol{\Sigma}_t^{1/2}-\overline{\boldsymbol{\Sigma}}^{1/2}{\mathfrak{O}}_{t})\mathbf{z}_t\mathbf{z}_t^T(\boldsymbol{\Sigma}_t^{1/2}-\overline{\boldsymbol{\Sigma}}^{1/2}{\mathfrak{O}}_{t})^T\Big)\Bigg)\\
=&\frac{1}{p}E\Bigg(\tr\Big(({\mathbf{W}}_t-\overline{\boldsymbol{\Sigma}}^{1/2})({\mathbf{W}}_t-\overline{\boldsymbol{\Sigma}}^{1/2})^T\Big)\Bigg)=o(1),
\endaligned
\end{equation}
which implies that
\begin{equation}\label{trace_diff_W}
\frac{1}{pn}\tr\Big(({{\mathbf{R}}}-\mathfrak{R}^{0})({{\mathbf{R}}}-\mathfrak{R}^{0})^T\Big)=o_p(1).
\end{equation}
The desired bound \eqref{L_Sw0_check} then follows from \eqref{L4_recovered_w}, \eqref{tr_tildS_w} and \eqref{trace_diff_W}. 

It remains to show  \eqref{E_tildR_w}.  
Define $$
\aligned
&\mathbf{y}_t=\overline{\boldsymbol{\Sigma}}^{-1/2}\mathbf{R}_t\text{, }\\
&\boldsymbol{\Sigma}_{y,t}=\overline{\boldsymbol{\Sigma}}^{-1/2}\boldsymbol{\Sigma}_t\overline{\boldsymbol{\Sigma}}^{-1/2}=(1-a-b)/(a-b)\mathbf{I}+\sum_{s=1}^\infty ab^{s-1} \mathbf{y}_{t-s}\mathbf{y}_{t-s}^T\\
&\mathbf{z}_t^y=(\boldsymbol{\Sigma}_{y,t})^{-1/2}\mathbf{y}_t=(\boldsymbol{\Sigma}_{y,t})^{-1/2}\overline{\boldsymbol{\Sigma}}^{-1/2}{\boldsymbol{\Sigma}}_t^{1/2}\mathbf{z}_t,\quad\text{ and}\\
&\mathcal{O}^y_t=(\boldsymbol{\Sigma}_{y,t})^{-1/2}\overline{\boldsymbol{\Sigma}}^{-1/2}{\boldsymbol{\Sigma}}_t^{1/2}.\endaligned$$
Note that $\mathcal{O}^y_t$ is $\mathcal{F}_{t-1}$-measurable, and $\mathcal{O}^y_t(\mathcal{O}^y_t)^T=(\mathcal{O}^y_t)^T\mathcal{O}^y_t=\mathbf{I}$. By the assumption that $\mathbf{z}_t\underset{\text{i.i.d.}}\sim N(0,\mathbf{I})$, we have that $\mathbf{z}_{t}^{y}$ is independent with $\boldsymbol{\Sigma}_{y,t}$, and
$\mathbf{z}_{t}^{y}\underset{\text{i.i.d.}}\sim N(0,\mathbf{I}).$
We then get that  $\mathbf{y}_t=(\boldsymbol{\Sigma}_{y,t})^{1/2}\mathbf{z}_t^y$ follows model \eqref{BEKK} with $E(\mathbf{y}_t\mathbf{y}_t^T)=E(\boldsymbol{\Sigma}_{y,t})=\mathbf{I}$. 
We have
\begin{equation}\label{WSgima}
 \aligned
 \tr(\mathbf{W}_t\overline{\boldsymbol{\Sigma}}^{1/2})=&\tr(\overline{\boldsymbol{\Sigma}}^{1/2}\boldsymbol{\Sigma}_{y,t}^{1/2}\overline{\boldsymbol{\Sigma}}^{1/2})\\
 =& \tr(\overline{\boldsymbol{\Sigma}})+\tr\big(\overline{\boldsymbol{\Sigma}}^{1/2}(\boldsymbol{\Sigma}_{y,t}^{1/2}-\mathbf{I})\overline{\boldsymbol{\Sigma}}^{1/2}\big).
 \endaligned
 \end{equation}
We denote by $x_1, ..., x_p$, and $u_1, ..., u_p$ the eigenvalues of $\boldsymbol{\Sigma}_{y,t}$ and the corresponding eigenvectors. By Assumption \ref{asump1}(iii) and the fact that $|\sqrt{x}-1|\leq |x-1|$ for all $x\geq 0$, we have
 $$
 \aligned
 &\frac{1}{p}\Bigg| \tr\Big(\overline{\boldsymbol{\Sigma}}^{1/2}(\boldsymbol{\Sigma}_{y,t}^{1/2}-\mathbf{I})\overline{\boldsymbol{\Sigma}}^{1/2}\Big)\Bigg|\\
 =&\frac{1}{p}\Bigg|\sum_{i=1}^p (\sqrt{x_i}-1) \Big(u_i^T\overline{\boldsymbol{\Sigma}}u_i\Big)\Bigg|\\
 \leq& C\sum_{i=1}^p |\sqrt{x_i}-1|/p\\
 \leq&  C\sqrt{\sum_{i=1}^p (x_i-1)^2/p}=C\sqrt{\tr(\boldsymbol{\Sigma}_{y,t}^2)/p-2\tr(\boldsymbol{\Sigma}_{y,t})/p+1}.
 \endaligned$$
 By Jensen's inequality and noting that $E(\mathbf{y}_t\mathbf{y}_t^T)=\mathbf{I}$, we get
\begin{equation}\label{iXt}
 \aligned
 &\frac{1}{p}E\Bigg(\Big| \tr\Big(\overline{\boldsymbol{\Sigma}}^{1/2}(\boldsymbol{\Sigma}_{y,t}^{1/2}-\mathbf{I})\overline{\boldsymbol{\Sigma}}^{1/2}\Big)\Big|\Bigg)\\
 \leq& C\sqrt{E\Big(\tr(\boldsymbol{\Sigma}_{y,t}^2)/p-2\tr(\boldsymbol{\Sigma}_{y,t})/p+1\Big)}\\
 = & C\sqrt{E\big(\tr(\boldsymbol{\Sigma}_{y,t}^2/p)\big)-1}.
 \endaligned
 \end{equation}
Next, we  compute $E\Big(\tr(\boldsymbol{\Sigma}_{y,t}^{2})\Big)$. 
Note that for any fixed matrix $\mathbf{A}\geq 0$ and  $\mathbf{z}\sim{N}(0, \mathbf{I})$,  we have
\begin{equation}\label{Eq2}
E\Big((\mathbf{z}^T\mathbf{A}\mathbf{z})^2\Big)=\big(\tr(\mathbf{A})\big)^2+2\tr(\mathbf{A}^2).
\end{equation}
Because $\boldsymbol{\Sigma}_{y,t+1}=(1-a-b)\mathbf{I}+a \mathbf{y}_{t}\mathbf{y}_{t}^T+ b\boldsymbol{\Sigma}_{y,t}$, we have 
\begin{equation}\label{trS2_y}
\aligned
\tr(\boldsymbol{\Sigma}_{y,t+1}^{2})=&(1-a-b)^2p+a^2 (\mathbf{y}_{t}^T\mathbf{y}_{t})^2+b^2\tr(\boldsymbol{\Sigma}_{y,t}^{2})\\
&+2(1-a-b) a \mathbf{y}_{t}^{T}\mathbf{y}_{t}\\
&+2(1-a-b) b \tr(\boldsymbol{\Sigma}_{y,t})\\
&+2a b \mathbf{y}_{t}^T\boldsymbol{\Sigma}_{y,t}\mathbf{y}_{t}.
\endaligned
\end{equation}
By normality of $\mathbf{z}^{y}_t$, the independence between $(\boldsymbol{\Sigma}_{y,t})$ and $(\mathbf{z}^y_t)$, \eqref{Eq2} and conditional on $\boldsymbol{\Sigma}_{y,t}$, we have 
$$
\aligned
E\Big((\mathbf{y}_{t}^T\mathbf{y}_{t})^2\Big)=&E\Big(\big((\mathbf{z}^y_{t})^T\boldsymbol{\Sigma}_{y,t}\mathbf{z}^y_{t}\big)^2\Big)\\
=&
E\Big(\big(\tr(\boldsymbol{\Sigma}_{y,t})\big)^2\Big)+2E\Big(\tr(\boldsymbol{\Sigma}^2_{y,t})\Big).
\endaligned
$$ Moreover, $E(\boldsymbol{\Sigma}_{y,t})=\mathbf{I}$. Taking expectations on both sides of \eqref{trS2_y} yields 
\begin{equation}\label{Eq1_S2_y}
\aligned
E\Big(\tr(\boldsymbol{\Sigma}_{y,t}^{2})\Big)=&\frac{1-(a+b)^2}{1-a^2-(a+b)^2}p\\
&+\frac{a^2}{1-a^2-(a+b)^2}E\Big(\big(\tr(\boldsymbol{\Sigma}_{y,t})\big)^2\Big).
\endaligned
\end{equation}

About $\Big(\tr(\boldsymbol{\Sigma}_{y,t})\Big)^2$, note that
\begin{equation}\label{Etr2_y}
\aligned
\Big(\tr(\boldsymbol{\Sigma}_{y,t+1})\Big)^{2}=&(1-a-b)^2 p^2+a^2 (\mathbf{y}_{t}^T\mathbf{y}_{t})^2+b^2\Big(\tr(\boldsymbol{\Sigma}_{y,t})\Big)^{2}\\
&+2(1-a-b) a p (\mathbf{y}_{t}^T\mathbf{y}_{t})\\
&+2(1-a-b) b p\tr(\boldsymbol{\Sigma}_{y,t})\\
&+2a b \mathbf{y}_{t}^T\mathbf{y}_{t}\tr(\boldsymbol{\Sigma}_{y,t}).
\endaligned
\end{equation}
Taking expectations on both sides of \eqref{Etr2_y} yields 
\begin{equation}\label{Eq2_S2_y}
\aligned
E\Bigg(\Big(\tr(\boldsymbol{\Sigma}_{y,t})\Big)^{2}\Bigg)=&p^2+\frac{2a^2}{1-(a+b)^2} E\Big(\tr(\boldsymbol{\Sigma}_{y,t}^2)\Big).
\endaligned
\end{equation}
Combining \eqref{Eq1_S2_y} and \eqref{Eq2_S2_y}, we get 
\begin{equation}\label{EM2_p3_y}
\aligned
&E\Big(\tr(\boldsymbol{\Sigma}_{y,t}^{2})\Big)=C_pC_p^{(2)}\Bigg(\frac{1-(a+b)^2}{a^2}p+p^2\Bigg),\text{ and}\\
&E\Bigg(\Big(\tr(\boldsymbol{\Sigma}_{y,t})\Big)^{2}\Bigg)=C_p\Bigg(p^2+2C_p^{(2)}p\Bigg),
\endaligned
\end{equation}
where
$$\aligned
&
C_p=\Bigg(1-\frac{2a^4}{\Big(1-a^2-(a+b)^2\Big)\Big(1-(a+b)^2\Big)}\Bigg)^{-1},\text{ and }\\
&C_p^{(2)}=\frac{a^2}{1-a^2-(a+b)^2},
\endaligned
$$
and \eqref{EM2_p3_y} holds when  $C_p$ and $C^{(2)}_p$ are both positive. Under Case II, 
we have 
$$\frac{1-(a+b)^2}{1-a^2-(a+b)^2}\to 1,\quad C_p=1+o(1),\quad \text{ and } \quad C_p^{(2)}= \frac{a^2}{1-(a+b)^2}\cdot\big(1+o(1)\big).
$$
It follows that
\begin{equation}\label{trX2}
E\big(\tr(\boldsymbol{\Sigma}_{y,t}^2)\big)=p+O\big(p^2a^2/(1-a-b) \big)=p+o(p).
\end{equation}
Combining \eqref{WSgima}, \eqref{iXt} and \eqref{trX2} yields
$$
E\Big(\tr(\mathbf{W}_t\overline{\boldsymbol{\Sigma}}^{1/2})\Big)=\tr(\overline{\boldsymbol{\Sigma}})+o(p).
$$
We then get
$$
\aligned
&\frac{1}{p}E\Bigg(\tr\Big(({\mathbf{W}}_t-\overline{\boldsymbol{\Sigma}}^{1/2})({\mathbf{W}}_t-\overline{\boldsymbol{\Sigma}}^{1/2})^T\Big)\Bigg)\\
=&\frac{1}{p}\Bigg(E\big(\tr({\mathbf{\Sigma}}_t)\big)-\tr(\overline{\boldsymbol{\Sigma}})\Bigg)+\frac{2}{p}\Bigg(\tr(\overline{\boldsymbol{\Sigma}})-E\Big(\tr({\mathbf{W}}_t\overline{\boldsymbol{\Sigma}}^{1/2})\Big)\Bigg)\\
=&o(1),
\endaligned
$$
namely, \eqref{E_tildR_w} holds. 
  $\qed$ 

\section*{Proof of~Theorem~\ref{nonreducable}}
Write 
$$
M^p_2=\frac{1}{p}\tr\big((\mathbf{S}_n)^2\big)=\frac{1}{n^2 p}\sum_{1\leq t_1, t_2\leq n}\mathbf{z}^T_{t_1}\boldsymbol{\Sigma}_{t_1}^{1/2} \boldsymbol{\Sigma}_{t_2}^{1/2}\mathbf{z}_{t_2}\mathbf{z}_{t_2}^T\boldsymbol{\Sigma}_{t_2}^{1/2}\boldsymbol{\Sigma}_{t_1}^{1/2}\mathbf{z}_{t_1}.
$$
We have
\begin{equation}\label{EM2}
\aligned
E(M^p_2)&
=\frac{1}{n^2 p}\sum_{t=1}^n E\Big((\mathbf{z}^T_{t}\boldsymbol{\Sigma}_{t}\mathbf{z}_{t})^2\Big) +  \frac{2}{n^2 p}\sum_{t_1=1}^{n-1}\sum_{t_2=t_1+1}^n E\Big(\mathbf{z}_{t_1}^T\boldsymbol{\Sigma}_{t_1}^{1/2}\boldsymbol{\Sigma}_{t_2}\boldsymbol{\Sigma}_{t_1}^{1/2}\mathbf{z}_{t_1}\Big)\\
&=:I+II.
\endaligned
\end{equation}
Under Assumption \ref{asump1}(i)--(iv), by (4.14) of \cite{yin1986limiting}, we have \begin{equation}\label{EM0}
E(M^{0,p}_2)=\tr\Big(\overline{\boldsymbol{\Sigma}}^2/p\Big)+y\Big(\tr(\overline{\boldsymbol{\Sigma}})/p\Big)^2+o(1).
\end{equation}
We will show
\begin{equation}\label{lowerbound_I_II}
\aligned
I&\geq y\Big(\tr(\overline{\boldsymbol{\Sigma}})/p\Big)^2\big(1+o(1)\big),\text{ and }II\geq \tr\Big(\overline{\boldsymbol{\Sigma}}^2/p\Big)\big(1+o(1)\big)+\delta. 
\endaligned 
\end{equation}
The desired bound \eqref{case1_thm2} then follows from \eqref{EM0} and \eqref{lowerbound_I_II}. 
 
We start with term $I$. By the independence between $(\boldsymbol{\Sigma}_t)$ and $(\mathbf{z}_t)$  and Jensen's inequality, we have
\begin{equation}\label{lowerbound_I}
\aligned I=&\frac{1}{n^2 p}\sum_{t=1}^n E\Big((\mathbf{z}^T_{t}\boldsymbol{\Sigma}_{t}\mathbf{z}_{t})^2\Big)\\
\geq &\frac{1}{n^2 p}\sum_{t=1}^n \Big(E(\mathbf{z}^T_{t}\boldsymbol{\Sigma}_{t}\mathbf{z}_{t})\Big)^2\\
=&\frac{1}{n p} \Big(\tr(\overline{\boldsymbol{\Sigma}})\Big)^2=y\Big(\tr(\overline{\boldsymbol{\Sigma}})/p\Big)^2\big(1+o(1)\big).
\endaligned
\end{equation}

About term $II$, by \eqref{Eq2}, the independence between $(\boldsymbol{\Sigma}_t)$ and $(\mathbf{z}_t)$, $\mathbf{z}_t\underset{\text{i.i.d.}}\sim N(0,\mathbf{I})$, taking conditional expectations recursively and using the fact that
$$
\aligned
\boldsymbol{\Sigma}_{t+j}=&\big(1-(a+b)^{j}\big)\overline{\boldsymbol{\Sigma}}\\
&+\sum_{i=1}^j a(a+b)^{j-i}(\boldsymbol{\Sigma}_{t+i-1}^{1/2}\mathbf{z}_{t+i-1}\mathbf{z}_{t+i-1}^T\boldsymbol{\Sigma}_{t+i-1}^{1/2}-\boldsymbol{\Sigma}_{t+i-1})\\
&+(a+b)^{j}\boldsymbol{\Sigma}_{t},
\endaligned
$$
we have
$$
\aligned
II=&\frac{2}{n^2 p}\sum_{t_1=1}^{n-1}\Big(\sum_{t_2=t_1+1}^n \big(1-(a+b)^{t_2-t_1}\big)\Big) E\Big(\mathbf{z}^T_{t_1}\boldsymbol{\Sigma}_{t_1}^{1/2}\overline{\boldsymbol{\Sigma}}\boldsymbol{\Sigma}_{t_1}^{1/2}\mathbf{z}_{t_1}\Big)\\
&+\frac{2a}{n^2 p}\sum_{t_1=1}^{n-1}\Big(\sum_{t_2=t_1+1}^n (a+b)^{t_2-t_1-1} \Big)\Bigg(E\Big((\mathbf{z}^T_{t_1}\boldsymbol{\Sigma}_{t_1}\mathbf{z}_{t_1})^2\Big)-E\Big(\mathbf{z}^T_{t_1}\boldsymbol{\Sigma}^2_{t_1}\mathbf{z}_{t_1}\Big)\Bigg)\\
&+\frac{2}{n^2 p}\sum_{t_1=1}^{n-1}\Big(\sum_{t_2=t_1+1}^n (a+b)^{t_2-t_1} \Big)E\Big(\mathbf{z}^T_{t_1}\boldsymbol{\Sigma}_{t_1}^{2}\mathbf{z}_{t_1}\Big)\\
=&\frac{2}{n^2 p}\sum_{t_1=1}^{n-1}\Big(\sum_{t_2=t_1+1}^n \big(1-(a+b)^{t_2-t_1}\big)\Big) \tr\big(\overline{\boldsymbol{\Sigma}}^2\big)\\
&+\frac{2}{n^2 p}\cdot\frac{a}{a+b}\cdot\sum_{t_1=1}^{n-1}\Big(\sum_{t_2=t_1+1}^n (a+b)^{t_2-t_1} \Big)E\Big(\big(\tr(\boldsymbol{\Sigma}_{t_1})\big)^{2}\Big)\\
&+\frac{2}{n^2 p}\cdot\Big(1+\frac{a}{a+b}\Big)\cdot\sum_{t_1=1}^{n-1}\Big(\sum_{t_2=t_1+1}^n (a+b)^{t_2-t_1} \Big)E\Big(\tr(\boldsymbol{\Sigma}_{t_1}^{2})\Big).
\endaligned$$
Define
$$
\aligned
\gamma_n=&\sum_{t_1=1}^{n-1}\Big(\sum_{t_2=t_1+1}^n (a+b)^{t_2-t_1}\Big)\\
=&\frac{(a+b)\Big((n-1)\big(1-(a+b)\big)-(a+b)+(a+b)^{n}\Big)}{\Big(1-(a+b)\Big)^2}.
\endaligned$$
It is easy to show that
\begin{equation}\label{vare_rate}
\gamma_n\asymp \frac{n}{1-a-b}\min\Big(n(1-a-b),1\Big).
\end{equation}
We then get 
\begin{equation}\label{EM2_p2}
\aligned
II=&\Big(\frac{n-1}{n}\Big)\tr\big(\overline{\boldsymbol{\Sigma}}^2\big)/p+\frac{2\gamma_n}{n^2p}\Bigg(E\Big(\tr(\boldsymbol{\Sigma}_{t}^{2})\Big)-\tr\big(\overline{\boldsymbol{\Sigma}}^2\big)\Bigg)\\
&+\frac{2 a\gamma_n }{(a+b)n^2p}E\Big(\big(\tr(\boldsymbol{\Sigma}_{t})\big)^{2}\Big)+\frac{2\gamma_n }{n^2p}\Big(\frac{a}{a+b}\Big) E\Big(\tr(\boldsymbol{\Sigma}_{t}^{2})\Big)\\
\geq &\frac{1}{p}\tr\big(\overline{\boldsymbol{\Sigma}}^2\big)\big(1+o(1)\big)+\frac{2\gamma_n}{n^2p}\Bigg(E\Big(\tr(\boldsymbol{\Sigma}_{t}^{2})\Big)-\tr\big(\overline{\boldsymbol{\Sigma}}^2\big)\Bigg).
\endaligned
\end{equation}
About $E\Big(\tr(\boldsymbol{\Sigma}_{t}^{2})\Big)$, by $\boldsymbol{\Sigma}_{t+1}=(1-a-b)\overline{\boldsymbol{\Sigma}}+a \boldsymbol{\Sigma}_{t}^{1/2} \mathbf{z}_{t}\mathbf{z}_{t}^T\boldsymbol{\Sigma}_{t}^{1/2}+ b\boldsymbol{\Sigma}_{t}$ and using an argument similar to \eqref{Eq2}, \eqref{trS2_y} and \eqref{Eq1_S2_y}, one can show that
when $a^2+(a+b)^2<1$,
\begin{equation}\label{Eq1_S2}
\aligned
E\Big(\tr(\boldsymbol{\Sigma}_{t}^{2})\Big)=&\frac{1-(a+b)^2}{1-a^2-(a+b)^2}\tr(\overline{\boldsymbol{\Sigma}}^2)\\
&+\frac{a^2}{1-a^2-(a+b)^2}E\Big(\big(\tr(\boldsymbol{\Sigma}_{t})\big)^2\Big).
\endaligned
\end{equation}
About $E\Big(\big(\tr(\boldsymbol{\Sigma}_{t})\big)^2\Big)$, by Jensen's inequality, we have
\begin{equation}\label{Eq2_S2}
\aligned
E\Big(\big(\tr(\boldsymbol{\Sigma}_{t})\big)^{2}\Big)\geq&\Big(E\big(\tr(\boldsymbol{\Sigma}_{t})\big)\Big)^{2}=\big(\tr(\overline{\boldsymbol{\Sigma}})\big)^2.
\endaligned
\end{equation}
Noting that $$\frac{a^2}{1-a^2-(a+b)^2}\geq \frac{a^2}{1-(a+b)^2}\geq \frac{a^2}{2(1-a-b)},$$ \eqref{Eq1_S2} and \eqref{Eq2_S2} imply that
\begin{equation}\label{EM2_p3_4}
\aligned
E\Big(\tr(\boldsymbol{\Sigma}_{t}^{2})\Big)\geq \tr(\overline{\boldsymbol{\Sigma}}^2)+\frac{a^2}{2(1-a-b)}\Big(\tr(\overline{\boldsymbol{\Sigma}})\Big)^2.
\endaligned
\end{equation}
By Assumption \ref{asump1}(ii), $\tr(\overline{\boldsymbol{\Sigma}})\asymp p$. By Assumption \ref{asump1}(iv), \eqref{vare_rate}, \eqref{EM2_p2} and \eqref{EM2_p3_4}, 
\begin{equation}\label{lowerbound_II}
II\geq  \tr\big(\overline{\boldsymbol{\Sigma}}^2/p\big)\big(1+o(1)\big)+C\Big(\eta(a,b,p)\Big)^2.
\end{equation}
The bound \eqref{lowerbound_I_II} follows from \eqref{lowerbound_I},  \eqref{lowerbound_II} and the assumption that  $\eta(a,b,p)>c$.
$\qed$

\section*{Proof of Corollary \ref{NLS_ADJ}}
Theorem \ref{recovered} implies that 
$F^{\widetilde{\mathbf{S}}_n}-F^{{\mathbf{S}}^0_n}\overset{\text{P}}\to 0$ as \mbox{$p, n\to \infty$}. By Theorem~1 of \cite{MP67}, 
$F^{{\mathbf{S}}^0_n}\overset{\text{P}}\to F$. 
Denote the empirical distribution of the truncated eigenvalues $(\widetilde{\lambda}_i^\tau)$ by $F^{\widetilde{\mathbf{S}}^\tau_n}$.  Then as long as the truncation level is greater than the upper bound of the support of $F$, we have
\begin{equation}\label{F_nconverge_F}
F^{\widetilde{\mathbf{S}}^\tau_n}\overset{\text{P}}\to F.
\end{equation}
The conclusion \hbox{$ \sum_{i=1}^p(\lambda^H_i-\widehat{\lambda}^H_i)^2/p=o_p(1)$} then follows from \eqref{F_nconverge_F} and the proof of~Theorem~2.2 of \cite{LW15}.~$\qed$

\section*{Proof of~Theorem~\ref{thetag}}\label{proof_thetag}
Note that
$$
\aligned
&\frac{1}{p}\tr\Big((\widetilde{\mathbf{S}}_n-z\mathbf{I})^{-1}g(\overline{\boldsymbol{\Sigma}})\Big)-\frac{1}{p}\tr\Big((\boldsymbol{\mathcal{S}}^{0}_n-z\mathbf{I})^{-1}g(\overline{\boldsymbol{\Sigma}})\Big)\\
=&\frac{1}{p}\tr\Big((\boldsymbol{\mathcal{S}}^{0}_n-z\mathbf{I})^{-1}(\boldsymbol{\mathcal{S}}^0_n-\widetilde{\mathbf{S}}_n)(\widetilde{\mathbf{S}}_n-z\mathbf{I})^{-1}g(\overline{\boldsymbol{\Sigma}})\Big),
\endaligned
$$
where, recall that $\boldsymbol{\mathcal{S}}^{0}_n=\boldsymbol{\mathcal{R}}^0(\boldsymbol{\mathcal{R}}^0)^T/n$ for $\boldsymbol{\mathcal{R}}^0$ defined in \eqref{R0*_def}. 
We have  $\|g(\overline{\boldsymbol{\Sigma}})\|<c$, $\|(\boldsymbol{\mathcal{S}}^{0}_n-z\mathbf{I})^{-1}\|<1/v$  and $\|(\widetilde{\mathbf{S}}_n-z\mathbf{I})^{-1}\|<1/v$, where $v$ is the imaginary part of $z$. By Corollary~A.12 and Theorem A.14 of \cite{BS10}, for some $C>0$,
\begin{equation}\label{sj}
\aligned
&\Bigg|\frac{1}{p}\tr\Big((\widetilde{\mathbf{S}}_n-z\mathbf{I})^{-1}g(\overline{\boldsymbol{\Sigma}})\Big)-\frac{1}{p}\tr\Big((\boldsymbol{\mathcal{S}}^{0}_n-z\mathbf{I})^{-1}g(\overline{\boldsymbol{\Sigma}})\Big)\Bigg|\\
\leq &\frac{1}{p}\sum_{i=1}^ps_i(\widetilde{\mathbf{S}}_n-\boldsymbol{\mathcal{S}}_n^{0})\cdot s_i\Big((\widetilde{\mathbf{S}}_n-z\mathbf{I})^{-1}g(\overline{\boldsymbol{\Sigma}})(\boldsymbol{\mathcal{S}}^{0}_n-z\mathbf{I})^{-1}\Big)\\
\leq &\frac{C}{p}\sum_{i=1}^p s_i(\widetilde{\mathbf{S}}_n-\boldsymbol{\mathcal{S}}_n^{0}),
\endaligned
\end{equation}
where for any matrix $\mathbf{A}$, $s_i(\mathbf{A})$ represents its $i$th singular value. 
Note that
$$
\aligned
\widetilde{\mathbf{S}}_n-\boldsymbol{\mathcal{S}}_n^{0}=\frac{1}{n}(\widetilde{\mathbf{R}}-\boldsymbol{\mathcal{R}}^0)(\widetilde{\mathbf{R}}-\boldsymbol{\mathcal{R}}^0)^T+\frac{1}{n}(\widetilde{\mathbf{R}}-\boldsymbol{\mathcal{R}}^0)(\boldsymbol{\mathcal{R}}^0)^T+\frac{1}{n}\boldsymbol{\mathcal{R}}^0(\widetilde{\mathbf{R}}-\boldsymbol{\mathcal{R}}^0)^T.
\endaligned
$$
By~Theorem~A.8 of \cite{BS10},
\begin{equation}\label{sj_0}
\aligned
\sum_{i=1}^p s_i(\widetilde{\mathbf{S}}_n-\boldsymbol{\mathcal{S}}_n^{0})
\leq& 2\sum_{i=1}^ps_i\Big(\frac{1}{n}(\widetilde{\mathbf{R}}-\boldsymbol{\mathcal{R}}^0)(\widetilde{\mathbf{R}}-\boldsymbol{\mathcal{R}}^0)^T\Big)\\
&+4\sum_{i=1}^ps_i\Big(\frac{1}{n}(\widetilde{\mathbf{R}}-\boldsymbol{\mathcal{R}}^0)(\boldsymbol{\mathcal{R}}^0)^T\Big).
\endaligned
\end{equation}
By \eqref{part2_check} and \eqref{trace_diff}, 
\begin{equation}\label{trace_diff_R0}
\aligned
&\frac{1}{pn}\tr\Big((\widetilde{\mathbf{R}}-\boldsymbol{\mathcal{R}}^0)(\widetilde{\mathbf{R}}-\boldsymbol{\mathcal{R}}^0)^T\Big)\\
\leq& \frac{2}{pn}\tr\Big((\check{{\mathbf{R}}}-\boldsymbol{\mathcal{R}}^0)(\check{{\mathbf{R}}}-\boldsymbol{\mathcal{R}}^0)^T\Big)+\frac{2}{pn}\tr\Big((\widetilde{\mathbf{R}}-\check{{\mathbf{R}}})(\widetilde{\mathbf{R}}-\check{{\mathbf{R}}})^T\Big)\\
=&o_p(1).
\endaligned
\end{equation}
Therefore,
\begin{equation}\label{sj_2}
\sum_{i=1}^ps_i\Big(\frac{1}{n}(\widetilde{\mathbf{R}}-\boldsymbol{\mathcal{R}}^0)(\widetilde{\mathbf{R}}-\boldsymbol{\mathcal{R}}^0)^T\Big)=\tr\Big((\widetilde{\mathbf{R}}-\boldsymbol{\mathcal{R}}^0)(\widetilde{\mathbf{R}}-\boldsymbol{\mathcal{R}}^0)^T/n\Big)=o_p(p).
\end{equation}
By the Cauchy-Schwarz inequality, we have
$$\aligned
\sum_{i=1}^p s_i\Big(\frac{1}{n}(\widetilde{\mathbf{R}}-\boldsymbol{\mathcal{R}}^0)(\boldsymbol{\mathcal{R}}^0)^T\Big)\leq& \sqrt{p}\sqrt{\sum_{i=1}^p\Bigg(s_i\Big(\frac{1}{n}(\widetilde{\mathbf{R}}-\boldsymbol{\mathcal{R}}^0)(\boldsymbol{\mathcal{R}}^0)^T\Big)\Bigg)^2}\\
=&\sqrt{p}\sqrt{\tr\Big(\frac{1}{n^2}(\widetilde{\mathbf{R}}-\boldsymbol{\mathcal{R}}^0)(\boldsymbol{\mathcal{R}}^0)^T\boldsymbol{\mathcal{R}}^0(\widetilde{\mathbf{R}}-\boldsymbol{\mathcal{R}}^0)^T\Big)}.
\endaligned
$$
By \eqref{invariant_normal} and Theorem 3.1 of \cite{yin1988limit}, $\|\sum_{t=1}^n\boldsymbol{\zeta}_t\boldsymbol{\zeta}_t^T/n\|=O_p(1)$. 
It follows that
\begin{equation}\label{l2_bound_s0}
\|(\boldsymbol{\mathcal{R}}^0)^T\boldsymbol{\mathcal{R}}^0/n\|=\|\boldsymbol{\mathcal{R}}^0(\boldsymbol{\mathcal{R}}^0)^T/n\|\leq \|\overline{\boldsymbol{\Sigma}}\|\cdot\Big\|\sum_{t=1}^n\boldsymbol{\zeta}_t\boldsymbol{\zeta}_t^T/n\Big\|=O_p(1).
\end{equation}
By \eqref{trace_diff_R0}, \eqref{l2_bound_s0},
 Corollary~A.12 and Theorem A.14 of \cite{BS10} again, 
\begin{equation}\label{sj_1}
\aligned
&\tr\Big(\frac{1}{n^2}(\widetilde{\mathbf{R}}-\boldsymbol{\mathcal{R}}^0)(\boldsymbol{\mathcal{R}}^0)^T\boldsymbol{\mathcal{R}}^0(\widetilde{\mathbf{R}}-\boldsymbol{\mathcal{R}}^0)^T\Big)\\
\leq& \|(\boldsymbol{\mathcal{R}}^0)^T\boldsymbol{\mathcal{R}}^0/n\|\cdot \tr\Big((\widetilde{\mathbf{R}}-\boldsymbol{\mathcal{R}}^0)(\widetilde{\mathbf{R}}-\boldsymbol{\mathcal{R}}^0)^T/n\Big)\\
=&o_p(p).
\endaligned
\end{equation}
Combining \eqref{sj}, \eqref{sj_0}, \eqref{sj_2} and \eqref{sj_1} yields
\begin{equation}\label{gESD_S0_Sn}
\Bigg|\frac{1}{p}\tr\Big((\widetilde{\mathbf{S}}_n-z\mathbf{I})^{-1}g(\overline{\boldsymbol{\Sigma}})\Big)-\frac{1}{p}\tr\Big((\boldsymbol{\mathcal{S}}^{0}_n-z\mathbf{I})^{-1}g(\overline{\boldsymbol{\Sigma}})\Big)\Bigg|=o_p(1).
\end{equation}
By Theorem 2 of \cite{LP11},   we have
\begin{equation}\label{gESD_S0}
\frac{1}{p}\tr\Big((\boldsymbol{\mathcal{S}}^{0}_n-z\mathbf{I})^{-1}g(\overline{\boldsymbol{\Sigma}})\Big)-\Theta^g(z)=o_p(1).
\end{equation}
The conclusion follows. 
  $\qed$
\section*{Proof of Theorem \ref{NLS_ADJ_cov}}
By Corollary \ref{NLS_ADJ}, the estimated population eigenvalues satisfy that $\sum_{i=1}^p(\widehat{\lambda}^H_i-\lambda^H_i)^2/p\overset{\text{P}}\to 0$. By Theorem~4 and Theorem~4 of \cite{LP11}, $\widetilde{\boldsymbol{\Sigma}}^{or}$ is the oracle shrinkage estimator. The desired result then follows from~the proof of~Theorem~3.1 of \cite{LW15}. $\qed$

\end{document}